\documentclass[journal,onecolumn,11pt]{IEEEtran}
\usepackage[cmex10]{amsmath}
\usepackage{amsthm}
\usepackage{amsmath}
\usepackage{amssymb}
\usepackage{amsfonts}
\usepackage{epic}
\usepackage{epsfig}
\usepackage{graphicx}
\usepackage[tight,footnotesize]{subfigure}
\usepackage{latexsym}
\usepackage{mathrsfs}
\usepackage{multirow}
\usepackage{newalg}
\usepackage[all]{xy}
\usepackage{cite}
\usepackage[usenames]{color}

\newtheorem{thm}{Theorem}
\newtheorem{prop}{Proposition}
\newcommand{\Amb}{\mathcal{A}}
\newcommand{\Alltop}{{\!\mbox{\footnotesize{A}}}}
\newcommand{\trans}{\mathrm{T}}
\newcommand{\ip}[2]{\langle #1, #2 \rangle}

\newcommand{\norm}[1]{\Vert #1 \Vert}

\newcommand{\normtwo}[1]{\Vert #1 \Vert_2}

\newcommand{\normone}[1]{\Vert #1 \Vert_1}

\newcommand{\abs}[1]{\vert #1 \vert}
\newcommand{\Abs}[1]{\Big\vert #1 \Big\vert}

\newcommand{\set}[1]{\{ #1 \}}

\newcommand{\gr}[1]{( #1 )}
\newcommand{\Gr}[1]{\Big( #1 \Big)}

\newcommand{\Real}{\mathbb R}
\newcommand{\C}{\mathbb C}

\newcommand{\eps}{\varepsilon}
\newcommand{\To}{\longrightarrow}

\newcommand{\bs}[1]{\boldsymbol{#1}}

\newcommand{\I}{\mathrm{i}} %%\newcommand{\I}{\mathfrak{i}}
\newcommand{\e}{\mathrm{e}} % exponential
\newcommand{\ith}{$i$th}
\newcommand{\jth}{$j$th}
\newcommand{\kth}{$k$th} %\newcommand{\kth}{$k^\textrm{th}$}

\newcommand{\Nth}{$N$th}
\newcommand{\half}{\frac{1}{2}}

\def\Pr{{\mathbb P}}

\newcommand{\cs}{compressed sensing}
\newcommand{\CS}{Compressed Sensing}
\newcommand{\Cs}{Compressed sensing}
\newcommand{\rep}{representation}
\newcommand{\up}{uncertainty principle}

\newcommand{\STFT}{short-time Fourier transform}
\newcommand{\tf}{time-frequency}

\newcommand{\TF}{Time-Frequency}
\newcommand{\ONB}{orthonormal basis}
\newcommand{\WLOG}{Without loss of generality}

\begin{document}
\bibliographystyle{IEEEtran} %\bibliographystyle{plain}
%%%%%%%%%%%%%%%%%%%%%%%%%%%%%%%%%%%%%%%%%%%%%%%%%%%%%%%%%%%%%%%%%%%%%%%%%%%%
%%%%%%%%%%%%%%%%%%%%%%%%%%%%%%%%%%%%%%%%%%%%%%%%%%%%%%%%%%%%%%%%%%%%%%%%%%%%
%%%%%%%%%%%%%%%%%%%%%%%%%%%%%%%%%%%%%%%%%%%%%%%%%%%%%%%%%%%%%%%%%%%%%%%%%%%%

\title{High-Resolution Radar\\ via Compressed Sensing}
\author{Matthew A. Herman and Thomas Strohmer
\thanks{The authors are with the Department of Mathematics, University of California,
Davis, CA 95616-8633, USA (\mbox{e-mail}: \texttt{\{mattyh, strohmer\}@math.ucdavis.edu}).}
\thanks{This work was partially supported by NSF Grant No. DMS-0511461 and NSF VIGRE Grant
Nos. DMS-0135345, DMS-0636297.}}

\maketitle

\begin{abstract}
A stylized compressed sensing radar is proposed in which the
time-frequency plane is discretized into an $N \times N$ grid.
Assuming the number of targets $K$ is small (i.e., $K \ll N^2$),
then we can transmit a sufficiently ``incoherent'' pulse and employ
the techniques of compressed sensing to reconstruct the target
scene. A theoretical upper bound on the sparsity $K$ is presented.
Numerical simulations verify that even better performance can be
achieved in practice. This novel compressed sensing approach offers
great potential for better resolution over classical radar.
\end{abstract}

\begin{IEEEkeywords}
Compressed sensing, radar, sparse recovery, matrix identification,
Gabor analysis, Alltop sequence.
\end{IEEEkeywords}

%%%%%%%%%%%%%%%%%%%%%%%%%%%%%%%%%%%%%%%%%%%%%%%%%%%%%%%%%%%%%%%%%%%%%%%%%%%%
%%%%%%%%%%%%%%%%%%%%%%%%%%%%%%%%%%%%%%%%%%%%%%%%%%%%%%%%%%%%%%%%%%%%%%%%%%%%
\section{Introduction}
%%%%%%%%%%%%%%%%%%%%%%%%%%%%%%%%%%%%%%%%%%%%%%%%%%%%%%%%%%%%%%%%%%%%%%%%%%%%
%%%%%%%%%%%%%%%%%%%%%%%%%%%%%%%%%%%%%%%%%%%%%%%%%%%%%%%%%%%%%%%%%%%%%%%%%%%%
%\subsection{Background and Motivation}
\IEEEPARstart{R}{adar}, sonar and similar imaging systems are in
high demand in many civilian, military, and biomedical applications.
The resolution of these systems is limited by classical \tf\ \up s.
Using the concepts of \cs, we propose a radically new approach to
radar, which under certain conditions provides better \tf\
resolution. In this simplified version of a monostatic,
single-pulse, far-field radar system we assume that the targets are
radially aligned with the transmitter and receiver. As such, we will
only be concerned with the range and velocity of the targets. Future
studies will include cross-range information.

There are three key points to be aware of with this approach:
(1)~The transmitted signal must be sufficiently ``incoherent.''
Although our results rely on the use of a deterministic signal (the
Alltop sequence), transmitting white noise would yield a similar
outcome. (2)~This approach does not use a matched filter. (3)~The
target scene is recovered by exploiting the imposed sparsity
constraints.

This report is a first step in formalizing the theory of \cs\ radar
and contains many assumptions. In particular, analog to digital
(A/D) conversion and related implementation details are ignored.
Some of these issues are discussed in \cite{Baraniuk} where the
potential to design simplified hardware is highlighted.

The rest of this section establishes notation and tools from \tf\
analysis, while Section~\ref{sect:CS} reviews the concepts of sparse
\rep s and \cs. Our main contribution can be found in
Sections~\ref{sect:Matrix_ID_via_CS}
and~\ref{sect:LTV_Application_Radar}. Other applications are
addressed in Section~\ref{sect:Other_Applications}.

%%%%%%%%%%%%%%%%%%%%%%%%%%%%%%%%%%%%%%%%%%%%%%%%%%%%%%%%%%%%%%%%%%%%%%%%%%%%
\subsection{Notation and Tools from \TF\ Analysis}
%%%%%%%%%%%%%%%%%%%%%%%%%%%%%%%%%%%%%%%%%%%%%%%%%%%%%%%%%%%%%%%%%%%%%%%%%%%%
In this paper boldface variables represent vectors and matrices,
while non-boldface variables represent functions with a
continuous domain. Throughout
this discussion we only consider %(possibly complex)
functions with finite energy, i.e., $f\in L^2(\Real)$. For two
functions $f,g\in L^2(\Real)$, their \emph{cross-ambiguity function}
of $\tau,\omega\in\Real$ is defined as~\cite{IMA}
\begin{equation} \label{eq:Cross_Ambiguity_Fnc}
\Amb_{fg}\gr{\tau,\omega} \;=\; \int_\Real f\gr{t + \tau/2}\:\!
\overline{g\gr{t - \tau/2}} \:\! \e^{-2\pi\I\omega t} dt,\\
\end{equation}
where $\overline{{\;\cdot\:}}$ denotes complex conjugation, and the
upright Roman letter $\I=\sqrt{-1}$. The
\emph{\STFT} (STFT) of $f$ with respect to $g$ is
$V_gf\gr{\tau,\omega} = \int_\Real f\gr{t}\: \overline{g\gr{t -
\tau}}\;\! \e^{-2\pi\I\omega t}dt.$ A simple change of variable
reveals that, within a complex factor, the cross-ambiguity function
is equivalent to the STFT
\begin{equation} \label{eq:Cross_Ambiguity_equiv_STFT}
\Amb_{fg}\gr{\tau,\omega} \;=\; \e^{\pi\I\omega\tau}\:\!
V_gf\gr{\tau,\omega}.
\end{equation}
When $f=g$ we have the \emph{(self) ambiguity function}
$\Amb_f\gr{\tau,\omega}$. The shape of the \emph{ambiguity surface}
$\abs{\Amb_f\gr{\tau,\omega}}$ of $f$ is bounded above the
\emph{\tf\ plane} $\gr{\tau,\omega}$ by
$\abs{\Amb_f\gr{\tau,\omega}} \leq\Amb_f\gr{0,0}=\normtwo{f}^2.$

The \emph{radar \up} \cite{Grochenig_AdvGabAnalysis} states that if
\begin{equation} \label{eq:Radar_Uncertainty_Princ}
\int\!\!\!\!\int_U \abs{\Amb_{fg}\gr{\tau,\omega}}^2 d\tau d\omega
\;\geq\; \gr{1-\varepsilon}\,\normtwo{f}^2\normtwo{g}^2
\end{equation}
for some \emph{support} $U\subseteq\Real^2$ and $\varepsilon\geq0$,
then the area
\begin{equation} \label{eq:Radar_Uncertainty_Princ_Conclusion}
\abs{U}\geq\gr{1-\varepsilon}.
\end{equation}
Informally, this can be interpreted as saying that the size of an
ambiguity function's ``footprint'' on the \tf\ plane can only be
made so small.

In classical radar, the ambiguity function of $f$ is the main factor
in determining the resolution between targets~\cite{Rihk}.
Therefore, the ability to identify two targets in the \tf\ plane is
limited by the essential support of $\Amb_f\gr{\tau,\omega}$ as
dictated by the radar uncertainty principle. The primary result of
this paper is that, under certain conditions, \cs\ radar achieves
better target resolution than classical radar.

%%%%%%%%%%%%%%%%%%%%%%%%%%%%%%%%%%%%%%%%%%%%%%%%%%%%%%%%%%%%%%%%%%%%%%%%%%%%
%%%%%%%%%%%%%%%%%%%%%%%%%%%%%%%%%%%%%%%%%%%%%%%%%%%%%%%%%%%%%%%%%%%%%%%%%%%%
\section{\CS} \label{sect:CS}
%%%%%%%%%%%%%%%%%%%%%%%%%%%%%%%%%%%%%%%%%%%%%%%%%%%%%%%%%%%%%%%%%%%%%%%%%%%%
%%%%%%%%%%%%%%%%%%%%%%%%%%%%%%%%%%%%%%%%%%%%%%%%%%%%%%%%%%%%%%%%%%%%%%%%%%%%
Recently, the signal processing/mathematics community has seen a
paradigmatic shift in the way information is represented, stored,
transmitted and recovered \cite{Don, CanRomTao, Tropp_Greed}. This
area is often referred to as \emph{Sparse Representations and \CS}.
Consider a discrete signal $\bs{s}$ of length $M$. We say that it is
\emph{$K$-sparse} if at most $K\ll M$ of its coefficients are
nonzero (perhaps under some appropriate change of basis). With this
point of view the \emph{true} information content of $\bs{s}$ lives
in at most $K$ dimensions rather than $M$. In terms of signal
acquisition it makes sense then that we should only have to measure
a signal $N\!\sim\!K$ times instead of $M$. We do this by making $N$
\emph{non-adaptive, linear observations} in the form of $\bs{y} =
\bs{\Phi s}$ where $\bs{\Phi}$ is a dictionary of size $N\times M$.
If $\bs{\Phi}$ is sufficiently ``incoherent,'' then the information
of $\bs{s}$ will be embedded in $\bs{y}$ such that it can be
perfectly recovered with high probability. Current reconstruction
methods include using greedy algorithms such as \emph{orthogonal
matching pursuit} (OMP) \cite{Tropp_Greed}, and solving the convex
problem:
\begin{equation} \label{eq:Basis_Pursuit}
\min{\norm{\bs{s'}}_1} \:\textrm{ s.t. }\: \bs{\Phi s'} = \bs{y}.
\end{equation}
The latter program is often referred to as \emph{Basis
Pursuit}\footnote{When in the presence of additive noise $\bs{e}$
the measurements are of the form $\bs{y} = \bs{\Phi s} + \bs{e}$. If
each element of the noise obeys $\abs{e_n} \le \eps$, then BP can be
reformulated as
$$
\min{\norm{\bs{s'}}_1} \:\textrm{ s.t. }\: \abs{\gr{\bs{\Phi s'} -
\bs{y}}_n} \leq \eps, \quad n=0,\ldots,N-1.$$}~(BP)~\cite{Don,
CanRomTao}. A new algorithm, \emph{regularized orthogonal matching
pursuit}~(ROMP)~\cite{NeedellVershynin} has recently been proposed
which combines the advantages of OMP with those of~BP.

%%%%%%%%%%%%%%%%%%%%%%%%%%%%%%%%%%%%%%%%%%%%%%%%%%%%%%%%%%%%%%%%%%%%%%%%%%%%
%%%%%%%%%%%%%%%%%%%%%%%%%%%%%%%%%%%%%%%%%%%%%%%%%%%%%%%%%%%%%%%%%%%%%%%%%%%%
\section{Matrix Identification via \CS} \label{sect:Matrix_ID_via_CS}
\subsection{Problem Formulation}
\begin{figure}[!htb]
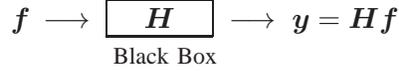

\begin{center}
\vspace{10pt}
    $\bs{f}\, \To \;\fbox{$\quad \bs{H} \quad $} \;\To \;\bs{y} = \bs{Hf}$\\
    \hspace{-34pt} {\footnotesize{Black Box}}
    \caption{Unknown system $\bs{H}$ with input probe
    $\bs{f}$ and output observation $\bs{y}$.} \label{fig:BlackBox}
\end{center}
\end{figure}

%%%%%%%%%%%%%%%%%%%%%%%%%%%%%%%%%%%%%%%%%%%%%%%%%%%%%%%%%%%%%%%%%%%%%%%%%%%%
%%%%%%%%%%%%%%%%%%%%%%%%%%%%%%%%%%%%%%%%%%%%%%%%%%%%%%%%%%%%%%%%%%%%%%%%%%%%
Consider an unknown matrix $\bs{H}\in\C^{N\times N'}$ and an \ONB\
(ONB) $\gr{\bs{H}_i}_{i=0}^{NN'\!-1}$ for $\C^{N\times N'}$. Note
that there are necessarily $NN'$ elements in this basis, and their
ortho-normality is with respect to the inner product derived from
the Frobenius norm (i.e., $\ip{\bs{A}}{\bs{B}}_F =
\textrm{trace}\gr{\bs{A}^*\bs{B}}$ for any
$\bs{A},\bs{B}\in\C^{N\times N'}$). Then there exist coefficients
$\gr{s_i}_{i=0}^{NN'\!-1}$ such that
\begin{equation} \label{eq:General_Unknown_Matrix_H}
\bs{H} \;=\, \sum_{i=0}^{NN'-1} s_i \bs{H}_i.
\end{equation}
Our goal is to identify/discover the coefficients
$\gr{s_i}_{i=0}^{NN'\!-1}$. Since the basis elements are fixed,
identifying these coefficients is tantamount to discovering
$\bs{H}$. We will do this by designing a test function $\bs{f} =
\gr{f_0,\dots,f_{N'-1}}^\trans\in~\C^{N'}$ and observing
$\bs{Hf}\in~\C^N$. Here, $\gr{\,\cdot\,}^\trans$ denotes the
transpose of a vector or a matrix. Figure~\ref{fig:BlackBox} depicts
this from a systems point of view where $\bs{H}$ is an unknown
``block box.'' Systems like this are ubiquitous in engineering and
the sciences. For instance, $\bs{H}$ may represent an unknown
communication channel which needs to be identified for equalization
purposes. In general, any linear time-varying (LTV) system can be
modeled by the basis of \tf\ shifts (described in the next section).

For simplicity, from now on assume that $N'=N$. The observation
vector can be reformulated as
\begin{equation} \label{eq:Reformulation_Observation}
\bs{y} \;=\; \sum_{i=0}^{N^2-1} s_i \bs{H}_i \bs{f} \;=\;
\sum_{i=0}^{N^2-1} s_i \bs{\varphi}_i \;=\; \bs{\Phi s},
\end{equation} where
\begin{equation} \label{eq:ith_Atom}
\bs{\varphi}_i = \bs{H}_i \bs{f} \in\C^N
\end{equation}
is the \emph{\ith\ atom}, $\bs{\Phi} = \gr{\bs{\varphi}_0 \,|\,
\cdots \,|\, \bs{\varphi}_{N^2-1}} \in\C^{N\times N^2}$ is the
concatenation of the atoms, and $\bs{s} = \gr{\:\!s_0, \cdots\!,
s_{N^2-1}\:\!}^\trans \in~\C^{N^2}$ is the coefficient vector. The
system of equations in (\ref{eq:Reformulation_Observation}) is
clearly highly underdetermined. If $\bs{s}$ is \emph{sufficiently
sparse}, then there is hope of recovering $\bs{s}$ from $\bs{y}$. To
use the reconstruction methods of \cs\ we need to design {$\bs{f}$}
so that the dictionary {$\bs{\Phi}$} is \emph{sufficiently
incoherent}.

%%%%%%%%%%%%%%%%%%%%%%%%%%%%%%%%%%%%%%%%%%%%%%%%%%%%%%%%%%%%%%%%%%%%%%%%%%%%
\subsection{The Coherence of a Dictionary}
%%%%%%%%%%%%%%%%%%%%%%%%%%%%%%%%%%%%%%%%%%%%%%%%%%%%%%%%%%%%%%%%%%%%%%%%%%%%
We are interested in how the atoms of a general dictionary
$\bs{\Phi} = \gr{\bs{\varphi}_i}_i \in\C^{N\times M}$ (with $N \leq
M$) are ``spread out'' in $\C^N$. This can be quantified by
examining the magnitude of the inner product between its atoms. The
\emph{coherence}~$\mu\gr{\bs{\Phi}}$ is defined as the maximum of
all of the distinct pairwise comparisons $\mu\gr{\bs{\Phi}} =
\max_{i\neq i'} \abs{\ip{\bs{\varphi}_i}{\bs{\varphi}_{i'}}}.$
Assuming that each $\normtwo{\bs{\varphi}_i} = 1$ the coherence is
bounded \cite{Ran55}, \cite{Wel74} by
\begin{equation} \label{eq:Coherence_Bound}
\sqrt{\frac{M-N}{N(M-1)}} \;\leq\; \mu(\bs{\Phi}) \;\leq\; 1.
\end{equation}
When $\mu(\bs{\Phi}) = 1$ we have two atoms which are aligned. This
is the worst-case scenario: \emph{maximal coherence}. In the other
extreme, when $\mu(\bs{\Phi}) = \sqrt{{(M-N)}/{N(M-1)}}$ we have
the best-case scenario: \emph{maximal incoherence}. Here the atoms
can be thought of as being ``spread out'' in $\C^N$. When a dictionary
can be expressed as the union of 2 or more ONBs,
this lower bound becomes ${1}/{\sqrt{N}}$ \cite{DGS75}.

%%%%%%%%%%%%%%%%%%%%%%%%%%%%%%%%%%%%%%%%%%%%%%%%%%%%%%%%%%%%%%%%%%%%%%%%%%%%
\subsection{The Basis of \TF\ Shifts}
%%%%%%%%%%%%%%%%%%%%%%%%%%%%%%%%%%%%%%%%%%%%%%%%%%%%%%%%%%%%%%%%%%%%%%%%%%%%
It is well-known from pseudo-differential operator theory
\cite{Grochenig} that any matrix can be represented by a basis of
\textbf{\tf\ shifts}. %(i.e., the composition of time shifts with
%frequency modulates).
Let the $N\times N$ matrices

$${\bs{T}} \;=\;
\left(
  \begin{array}{cccc}
    \!0  &          &          &  \!1 \! \\
    \!1  &     0    &          &       \\
         &  \ddots  &  \!\ddots  &       \\
    \!0  &          &    \! 1    &  \!0\!
  \end{array}
\right),
\quad
    {\bs{M}} \;=\;
\left(
  \begin{array}{cccc}
    \!\omega_N^0  &                   &              & \!\!\!\!0       \\
                    &  \!\!\!\!\omega_N^1   &              &                 \\
                    &                   &  \!\!\!\!\ddots  &                 \\
    \!\!0           &                   &              &  \!\!\!\!\omega_N^{N\!-\!1\!\!\!\!}
  \end{array}
\right)
$$

\vspace{10pt} \noindent respectively denote the \emph{unit-shift}
and \emph{modulation operators} where $\omega_N = \e^{2\pi\I/N}$ is
the \Nth\ root of unity. The \ith\ \tf\ basis element is defined as
\begin{equation} \label{eq:TF_Basis_Element_H_i}
\bs{H}_i \;=\; \bs{M}^{i \!\!\!\!\!\mod{\! N}} \!\cdot
\bs{T}^{\lfloor i/N\rfloor},
\end{equation}
where $\lfloor\cdot\rfloor$ is the floor function. A simple
calculation shows that the family $\gr{\bs{H}_i}_{i=0}^{N^2-1}$
forms an ONB with respect to the Frobenius inner product. Further,
under this basis it is known that some practical systems $\bs{H}$
with meaningful applications have a sparse \rep\
$\bs{s}$~\cite{Bel63, LiPreisig, ShanSwaryAviyente}. This fact
complements the theorems developed in the subsequent sections.

A finite collection of length-$N$ vectors which are \tf\ shifts of a
generating vector, and which spans the space $\C^N$ is called a
(discrete) \emph{Gabor frame} \cite{Grochenig}.
Since~$\gr{\bs{H}_i}_{i=0}^{N^2-1}$ is an ONB, it follows that our
dictionary~$\bs{\Phi}$ is a Gabor frame. \WLOG, assume
$\normtwo{\bs{f}}=1$. Because each $\bs{H}_i$ is a unitary matrix we
have from~(\ref{eq:ith_Atom}) that $\normtwo{\bs{\varphi}_i}=1$ for
$i=0,\ldots, N^2-1$. We can also express $\bs{\Phi}$ as the
concatenation of $N$ blocks
\begin{equation} \label{eq:Gabor_Dictionary_N_blocks}
\bs{\Phi} \;=\; \Gr{\,\bs{\Phi}^{(0)} \,|\, \bs{\Phi}^{(1)} \,|\,
\cdots \,|\, \bs{\Phi}^{(N-1)}\,},
\end{equation}
where the \kth\ block $\bs{\Phi}^{(k)} =
\bs{D}_{k}\!\cdot\!\bs{W}_{\!\!N}$, with $\bs{D}_{k}=
\mbox{diag}\set{f_k, \ldots, f_{N-1}, f_0, \ldots, f_{k-1}}$, and
$\bs{W}_{\!\!N} = \gr{\omega_N^{pq}}_{p,q=0}^{N-1}$. Here,
$\bs{\Phi}^{(k)}, \bs{D}_{k},$ and $\bs{W}_{\!\!N}$ are all matrices
of size $N\times N$. Essentially, the first column of
$\bs{\Phi}^{(k)}$ consists of the vector $\bs{f}$ shifted by $k$
units in time (with no modulation). The remaining $N\!-\!1$ columns
of $\bs{\Phi}^{(k)}$ consist of the $N\!-\!1$ other possible
modulations of this first column. Since there are $N$ different
modulates for each of the $N$ time shifts, we have $N^2$
combinations of \tf\ shifts, and these form the atoms of our
dictionary.

%%%%%%%%%%%%%%%%%%%%%%%%%%%%%%%%%%%%%%%%%%%%%%%%%%%%%%%%%%%%%%%%%%%%%%%%%%%%
\subsection{The Probing Test Function $\bs{f}$}
%%%%%%%%%%%%%%%%%%%%%%%%%%%%%%%%%%%%%%%%%%%%%%%%%%%%%%%%%%%%%%%%%%%%%%%%%%%%
We now introduce a candidate probe function $\bs{f}$ which results
in remarkable incoherence properties for the dictionary~$\bs{\Phi}$.
Consider the \emph{Alltop sequence} $\bs{f}_\Alltop =
\gr{f_n}_{n=0}^{N-1}$ for some prime $N\geq5$, where \cite{Alltop}
\begin{equation} \label{eq:Alltop_sequnce}
f_n \;=\; \frac{1}{\sqrt{N}}\:\!\e^{2\pi\I n^3/N}.
\end{equation} This
function has been proposed for use in telecommunications (CDMA,
etc.), for constructing the \emph{mutually unbiased bases} (MUBs)
used in quantum physics and quantum cryptography~\cite{Quantum_MUB},
and was made popular in the frames community
in~\cite{StrohmerHeath}.

Let $\bs{\Phi}_\Alltop$ denote the Gabor frame generated by the
Alltop sequence (\ref{eq:Alltop_sequnce}). Since its atoms are
already grouped into $N\times N$ blocks in
(\ref{eq:Gabor_Dictionary_N_blocks}), we will maintain this
structure by denoting the \jth\ atom of the \kth\ block as
$\bs{\varphi}_j^{(k)}$. %for any $j,k =~0, \ldots, N-1$.
Note that $\normtwo{\bs{f}_\Alltop} = 1$, so we have
$0 \leq \abs{\ip{\bs{\varphi}^{(k)}_j}{\bs{\varphi}^{(k')}_{j'}}} \leq 1$
for any $j,j',k,k' = 0, \ldots, N-1$. Within the \emph{same} block
(i.e., $k=k'$) we have
$$ \mbox{\textbf{Property 1:}} \qquad \abs{\ip{\bs{\varphi}^{(k)}_j}{\bs{\varphi}^{(k)}_{j'}}}
\;=\;
    \left\{
    \begin{array}{cl}
        0, & \mbox{  if $j\neq j'$}\\
        1, & \mbox{  if $j = j'$.}
    \end{array}
    \right. \label{property:One}
$$
Thus, each $\bs{\Phi}^{(k)}$ is an ONB for $\C^N$. Moreover, for
\emph{different} blocks (i.e., $k\neq k'$) we have
$$ \hspace{-55pt} \mbox{\textbf{Property 2:}}
\qquad \abs{\ip{\bs{\varphi}^{(k)}_j}{\bs{\varphi}^{(k')}_{j'}}}
\;=\; \displaystyle \frac{1}{\sqrt{N}}$$ \textbf{for all}
$j,j'=0,\ldots,N-1$. This means that there is a \emph{mutual
incoherence} between the atoms of different blocks (equivalently,
the $N$ blocks make up a set of MUBs). Trivially, it follows that
$\mu(\bs{\Phi}_\Alltop) = {1}/{\sqrt{N}}$. Furthermore, with
$M=N^2$ in (\ref{eq:Coherence_Bound}) we see that the lower bound of
${1}/{\sqrt{N+1}}$ is \emph{practically attained}.
These amazing properties are due to the cubic phase factor in the
Alltop sequence (\ref{eq:Alltop_sequnce}), and the fact that $N$ is
prime. More details and proofs can be found in \cite{Alltop}.

\emph{\textbf{Remark}}. Actually, in theory the Alltop sequence
yields a set of $N+1$ MUBs. This can be achieved by adjoining the
$N$ canonical unit vectors to the $N^2$ time-frequency shifted
Alltop sequences. This results in a total of $N^2+N$ vectors
(grouped in $N+1$ MUBs) that still maintain Properties~1 and~2.
However, this last MUB is simply the identity matrix. Since it
possesses no intrinsic \tf\ structure, we do not see how to use this
fact to our advantage in the context of radar.

\emph{\textbf{Remark}}. By inspection of (\ref{eq:Coherence_Bound})
we observe that the smallest possible incoherence for $M=N^2$
vectors is ${1}/{\sqrt{N+1}}$ which is slightly smaller than the
incoherence of the Gabor frame resulting from the Alltop sequence.
If a set of vectors obtains this optimal bound, it is automatically
an equiangular tight frame, see~\cite{StrohmerHeath}. It is
conjectured that for any $N$ there exists an (equiangular tight)
Gabor frame with $N^2$ elements which achieves the bound
${1}/{\sqrt{N+1}}$. However, explicit constructions are known only
for a very few cases, cf.~\cite{Appleby}. Therefore, and because the
difference between ${1}/{\sqrt{N}}$ and ${1}/{\sqrt{N+1}}$ is
negligible for large $N$, we will continue our investigation using
Alltop sequences.

%%%%%%%%%%%%%%%%%%%%%%%%%%%%%%%%%%%%%%%%%%%%%%%%%%%%%%%%%%%%%%%%%%%%%%%%%%%%
\subsection{Identifying Matrices via \CS: Theory}
%%%%%%%%%%%%%%%%%%%%%%%%%%%%%%%%%%%%%%%%%%%%%%%%%%%%%%%%%%%%%%%%%%%%%%%%%%%%
Having established the incoherence properties of the dictionary
$\bs{\Phi}_\Alltop$ we can now move on to apply the concepts and
techniques of \cs. It is worth pointing out that most \cs\ scenarios
deal with a $K$-sparse signal $\bs{s}$ (for some \emph{fixed} $K$),
and one is tasked with determining how many observations are
necessary to recover the signal. Our situation is markedly
different. Due to the fact that $\bs{\Phi}_\Alltop$ is constrained
to be $N\times N^2$, we know $\bs{y} = \bs{\Phi}_\Alltop\bs{s}$ will
contain exactly $N$ observations. With $N$ fixed, our \cs\ dilemma
is to determine how sparse $\bs{s}$ should be such that it can be
recovered from $\bs{y}$.

Therefore, with $N$ measurements, we can only consider recovering
signals which are less than $N$-sparse. Indeed, we hope to recover
any $K$-sparse signal $\bs{s}$ with $K\leq C\cdot N/\log{N}$ for
some $C>0$. The following theorems summarize the recovery of
$N\times N$ matrices via \cs\ when identified with the Alltop
sequence. Their proofs appear in Appendix~\ref{app:Proofs}. Assume
throughout that prime~$N\geq5$.

\begin{thm} \label{thm:Matrix_Recovery_Guaranteed}
Suppose $\bs{H} = \sum_i s_i \bs{H}_i \in\C^{N\times N}$ has a
$K$-sparse \rep\ under the \tf\ ONB, with $K < \half
\gr{\sqrt{N}+1}$, and that we have observed $\bs{y} =
\bs{H}\bs{f}_\Alltop$. Then we are \textbf{guaranteed} to
recover~$\bs{s}$ either via BP or OMP.
\end{thm}

The sparsity condition in
Theorem~\ref{thm:Matrix_Recovery_Guaranteed} is rather strict.
Instead of the requirement of \emph{guaranteed} perfect recovery, we
can ask to achieve it with only \emph{high probability}. This more
modest expectation provides us with a sparsity condition which is
more generous.

Unless specified otherwise, a \emph{\textbf{random}} signal in this
paper refers to a \textbf{vector} whose nonzero (complex)
coefficients are independent with a Gaussian distribution of zero
mean and unit variance.\footnote{For complex signals, each nonzero
entry has real and imaginary parts which are independent, Gaussian
random variables with zero mean and a variance of $1/2$; thus the
unit variance of each nonzero coefficient is the result of the sum
of the variances of its real and imaginary parts. From the
rotational invariance of the Gaussian distribution it can be shown
that the phase of each random coefficient is circularly symmetric,
i.e., its phase is uniformly distributed on the interval $[0,2\pi)$.
See Appendix~A of \cite{TseViswanath}.} Further, these nonzero
coefficients are uniformly distributed along the length of the
vector.

\begin{thm} \label{thm:Matrix_Recovery_High_Prob}
Suppose random $\bs{s}\in \C^{N^2}$ is a $K$-sparse vector with $K
\leq N/\gr{16\log{\gr{N/\varepsilon}}}$ for some sufficiently
small~$\varepsilon$. Suppose further that $\bs{H} = \sum_i s_i
\bs{H}_i \in \C^{N\times N}$ and that we have observed $\bs{y}
=\bs{H}\bs{f}_\Alltop$. Then BP will recover~$\bs{s}$ with
probability greater than $1-2\varepsilon^2-K^{-\vartheta}$ for some
$\vartheta\geq1$ s.t. $\sqrt{\vartheta
\log{N}/\log{\gr{N/\varepsilon}}} \leq c$ where $c$ is an absolute
constant.
\end{thm}

%%%%%%%%%%%%%%%%%%%%%%%%%%%%%%%%%%%%%%%%%%%%%%%%%%%%%%%%%%%%%%%%%%%%%%%%%%%%
\emph{\textbf{With Additive Noise.}}
%%%%%%%%%%%%%%%%%%%%%%%%%%%%%%%%%%%%%%%%%%%%%%%%%%%%%%%%%%%%%%%%%%%%%%%%%%%%
Theorems~\ref{thm:Matrix_Recovery_Guaranteed} and
\ref{thm:Matrix_Recovery_High_Prob} can be extended to include the
case of noisy observed signals. This will of course have an effect
on the sparsity of the signal of interest. For instance, the value
of $K$ in Theorem~\ref{thm:Matrix_Recovery_Guaranteed} is reduced
from $\half \gr{\sqrt{N}+1}$ to $\half \gr{\sqrt{N}+1}/\gr{1 + 2\eps
N/T}$ as seen in the following theorem.

\begin{thm}
\label{thm:Matrix_Recovery_Guaranteed_Noise} Suppose $\bs{H} =
\sum_i s_i \bs{H}_i \in\C^{N\times N}$ has a $K$-sparse \rep\ under
the \tf\ ONB, with $K < \half \gr{\sqrt{N}+1}/\gr{1 + 2\eps N/T}$.
Suppose further that we have observed $\bs{y} = \bs{H}\bs{f}_\Alltop
+ \bs{e}$, where each element of the noise $\abs{e_n} \le \eps$.
Then the solution $\bs{s^\bigstar}$ to BP exhibits stability
$\normone{\bs{s}-\bs{s^\bigstar}}\le T$.
\end{thm}

In a similar way, Theorem~\ref{thm:Matrix_Recovery_High_Prob} can be
rephrased to account for observed signals which have been perturbed.

%%%%%%%%%%%%%%%%%%%%%%%%%%%%%%%%%%%%%%%%%%%%%%%%%%%%%%%%%%%%%%%%%%%%%%%%%%%%
\subsection{Identifying Matrices via \CS: Simulation}
%%%%%%%%%%%%%%%%%%%%%%%%%%%%%%%%%%%%%%%%%%%%%%%%%%%%%%%%%%%%%%%%%%%%%%%%%%%%
Numerical simulations were performed and indicate that the theories
above are actually somewhat pessimistic. The simulations were
conducted as follows. The values of prime~$N$ ranged from $5$ to
$127$, and the sparsity $K$ ranged from $1$ to~$N$. For each ordered
pair $\gr{N,K}$ a complex-valued, $K$-sparse vector $\bs{s}$ of
length~$N^2$ was randomly generated. With this random signal the
observation $\bs{y}=\bs{\Phi}_\Alltop\bs{s}$ was generated. Then,
$\bs{y}$ and $\bs{\Phi}_\Alltop$ were input to convex optimization
software~\cite{Boyd, Mosek} to implement
BP~(\ref{eq:Basis_Pursuit}). Denote $\bs{s^\bigstar}$ as the
solution to the BP program. The recovered vector was deemed
successful if the error $\normtwo{\bs{s}-\bs{s^\bigstar}} \leq
10^{-4}$.  This procedure was repeated $100$ times for each
$\gr{N,K}$-pair; the total number of successes was recorded and then
averaged.

\begin{figure}[!ht]
\begin{center}
\includegraphics[width=3.5in]{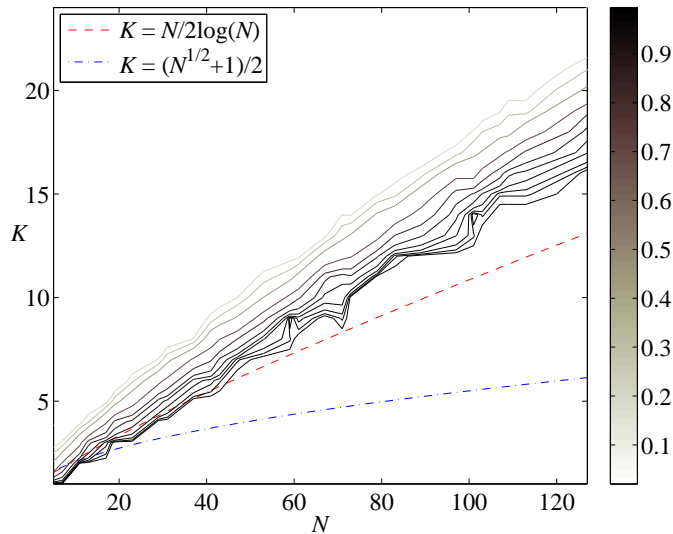}
    \caption{Numerical results from 100 independent Matlab simulations implementing
    BP for different $\gr{N,K}$-pairs.
    The solid, gray-black lines are contours whose values represent
    the fraction of successful recoveries
    vs. the $N$-$K$ domain. The dashed, red line shows that Theorem
    \ref{thm:Matrix_Recovery_High_Prob}
    is overly pessimistic. The region below this is the zone of
    ``perfect reconstruction.'' The lower dash-dotted, blue line illustrates
    that Theorem~\ref{thm:Matrix_Recovery_Guaranteed} is too strict.}
    \label{fig:AlltopSurfAndContourPhasePlot}
\end{center}
\end{figure}

Figure~\ref{fig:AlltopSurfAndContourPhasePlot} shows how the
numerical simulations compare to
Theorems~\ref{thm:Matrix_Recovery_Guaranteed} and
\ref{thm:Matrix_Recovery_High_Prob}. The fraction of successful BP
recoveries as a function of $\gr{N,K}$ is shown as solid, gray-black
contour lines. Although the values of $N$ used in the simulations
were relatively small, we see from these numerical results what
appears to be a trend. The dashed, red line represents $K =
N/\gr{2\log{N}}$, and the zone of ``perfect
reconstruction'' lies \emph{below} this line. %\footnote{ In previous
%simulations with \emph{real} values of $\bs{s}$ we found that the
%zone of ``perfect reconstruction'' was bounded by $K = N/\log{N}$.}
In this region a random $N\times N$ matrix (i.e., $\bs{H}$ as
defined in Theorem~\ref{thm:Matrix_Recovery_High_Prob}) with $1 \leq
K \leq N/\gr{2\log{N}}$ can be perfectly recovered with high
probability by observing $\bs{y} =\bs{H}\bs{f}_\Alltop$. This is
empirical evidence that the denominator in the upper bound of $K$ in
Theorem~\ref{thm:Matrix_Recovery_High_Prob} can be relaxed from
$\log{\gr{N/\varepsilon}}$ to just $\log{N}$, and that the
proportionality constant $C=1/2$. However, it is still an open
mathematical problem to prove this for the Alltop sequence.
Furthermore, the overly strict constraint of
Theorem~\ref{thm:Matrix_Recovery_Guaranteed} can be seen by the
lower dash-dotted, blue line representing
$K~=~\half\gr{\sqrt{N}+1}$.

%%%%%%%%%%%%%%%%%%%%%%%%%%%%%%%%%%%%%%%%%%%%%%%%%%%%%%%%%%%%%%%%%%%%%%%%%%%%
%%%%%%%%%%%%%%%%%%%%%%%%%%%%%%%%%%%%%%%%%%%%%%%%%%%%%%%%%%%%%%%%%%%%%%%%%%%%
\section{Radar} \label{sect:LTV_Application_Radar}
\subsection{Classical Radar Primer} \label{sect:Classical_Radar}
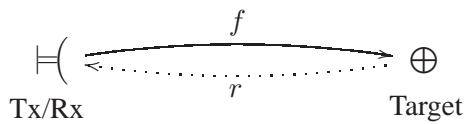
\begin{figure}[!htb]
\begin{center}
\begin{displaymath}
    \xymatrix{
        \models\!\!\!\Big(
        \;\ar@/^/@{->}[rrrr]^{\mbox{${f}$}} & & & &
        \ar@/^/@{.>}[llll]^{\!\mbox{${r}$}}
        \;\bigoplus
    }
    \vspace{-15pt}
\end{displaymath}
\begin{displaymath}
    \xymatrix{
        \mbox{Tx/Rx}& & & &
        \quad \!\!\!\!\!\!\!\!
        \!\!\mbox{Target}
    }
\end{displaymath}
    \caption{Simplified radar model. Tx transmits signal~$f$, and
    Rx receives the reflected (or echoed) signal~$r$ according to
    (\ref{eq:Radar_reflected_signal}).} \label{fig:Simple_Radar}
\end{center}
\end{figure}

%%%%%%%%%%%%%%%%%%%%%%%%%%%%%%%%%%%%%%%%%%%%%%%%%%%%%%%%%%%%%%%%%%%%%%%%%%%%
%%%%%%%%%%%%%%%%%%%%%%%%%%%%%%%%%%%%%%%%%%%%%%%%%%%%%%%%%%%%%%%%%%%%%%%%%%%%
Consider the following simple (narrowband) 1-dimensional,
monostatic, single-pulse, far-field radar model. \emph{Monostatic}
refers to the setup where the transmitter (Tx) and receiver (Rx) are
collocated. The far-field assumption permits us to model the targets
as point sources.
%\footnote{Note that the analysis presented in this paper can be
%easily extended to multiple Tx, Rx, which are moving as well as to
%pulse-train radar systems}
Suppose a target located at \emph{range} $x$ is traveling with
\emph{constant velocity} $v$ and has \emph{reflection coefficient}
$s_{xv}$. Figure~\ref{fig:Simple_Radar} shows such a radar with one
target. After transmitting signal $f(t)$, the receiver observes the
reflected signal
\begin{equation} \label{eq:Radar_reflected_signal}
r(t) \;=\; s_{xv}\:\! f(t-\tau_x)\:\!\e^{2\pi\I\omega_vt},
\end{equation}
where $\tau_x = 2x/c$ is the round trip time of flight, $c$ is the
speed of light, $\omega_v \approx -2\omega_0 v/c$ is the Doppler
shift, and~$\omega_0$ is the carrier frequency. The basic idea is
that the \emph{range-velocity} information $\gr{x,v}$ of the target
can be inferred from the observed \emph{time delay-Doppler shift}
$\gr{\tau_x,\omega_v}$ of $f$ in (\ref{eq:Radar_reflected_signal}).
Hence, a \tf\ shift operator basis is a natural \rep\ for radar
systems~\cite{AuslanderTolimieri}.

%\vspace{15pt}

\begin{figure}[!ht]
\begin{center}
\hspace{-125pt} \setlength{\unitlength}{4pt}
\begin{picture}(0,35)(0,0)
    \multiput(0,0)(2,0){18} {\line(0,1){34}}
    \multiput(0,0)(0,2){18} {\line(1,0){34}}
    \put(-2,-2.5){{$0$}}
    \put(31,-2.75){{$N\!-\!1$}}
    \put(14,-4){$\tau \rightarrow$}
    \put(-6.25,33.25){{$N\!-\!1$}}
    \put(-3.5,15){$\omega$}
    \put(-3.35,18){$\uparrow$}
    \put(9.4,17.4){$\bullet$}
    \put(21.4,11.4){$\bullet$}
    \put(19.4,9.4){$\bullet$}
    \put(17.4,23.4){$\bullet$}
    \put(25.4,19.4){$\bullet$}
    \put(10,18){\circle{7.45}}
    \put(22,12){\circle{7.45}}
    \put(20,10){\circle{7.45}}
    \put(18,24){\circle{7.45}}
    \put(26,20){\circle{7.45}}
\end{picture}
\vspace{20pt}
    \caption{The \tf\ plane discretized into an $N\times N$
    grid. Shown are five targets with their associated uncertainty regions.
    Classical radar detection techniques may fail to
    resolve the two targets whose regions are intersecting.
    In contrast, \cs\ radar will be able to distinguish them as long as
    the total number of targets is much less then $N^2$.}
    \label{fig:CS_Time_Freq_Plane}
\end{center}
\end{figure}
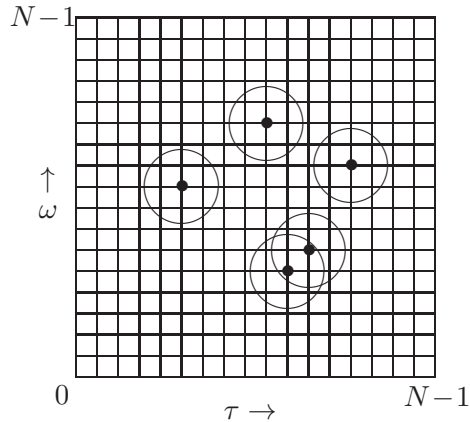

Using a \textbf{matched filter} at the receiver, the reflected
signal~$r$ is correlated with a \tf\ shifted version of the
transmitted signal~$f$ via the cross-ambiguity
function~(\ref{eq:Cross_Ambiguity_Fnc})
\begin{eqnarray} \label{eq:Matched_Filter_Ambiguity_Eqn}
\abs{\Amb_{rf}\gr{\tau,\omega}} &=& \Abs{\! \int_\Real
\!r(t)\:\!\overline{f(t-\tau)}\:\!\e^{-2\pi\I\omega t} dt\:\!} \nonumber\\
    &=& \abs{\:\! s_{xv} \:\! V_f\!f\gr{\tau-\tau_x, \omega-\omega_v}} \nonumber\\
    &=& \abs{\:\! s_{xv} \:\! \Amb_f\gr{\tau-\tau_x, \omega-\omega_v}}.
\end{eqnarray}
From this we see that the \tf\ plane consists of the ambiguity
surface of $f$ \emph{centered at the target's ``location''}
$\gr{\tau_x,\omega_v}$ and scaled by its reflection coefficient
$\abs{s_{xv}}$. Extending (\ref{eq:Matched_Filter_Ambiguity_Eqn}) to
include multiple targets is straightforward.
Figure~\ref{fig:CS_Time_Freq_Plane} illustrates an example of the
\tf\ plane with five targets; two of these have overlapping
uncertainty regions. The uncertainty region is a rough indication of
the essential support of $\Amb_f$
in~(\ref{eq:Radar_Uncertainty_Princ}). Targets which are too close
will have overlapping ambiguity functions. This may blur the exact
location of a target, or make uncertain how many targets are located
in a given region in the \tf\ plane. Thus, the range-velocity
resolution between targets of classical radar is limited by the
radar \up.

%%%%%%%%%%%%%%%%%%%%%%%%%%%%%%%%%%%%%%%%%%%%%%%%%%%%%%%%%%%%%%%%%%%%%%%%%%%%
\subsection{\CS\ Radar} \label{sect:CS_Radar}
%%%%%%%%%%%%%%%%%%%%%%%%%%%%%%%%%%%%%%%%%%%%%%%%%%%%%%%%%%%%%%%%%%%%%%%%%%%%
We now propose our stylized \cs\ radar which under appropriate
conditions can ``beat'' the classical radar \up! Consider $K$
targets with unknown range-velocities and corresponding reflection
coefficients. Next, discretize the \tf\ plane into an $N\times N$
grid as depicted in Figure~\ref{fig:CS_Time_Freq_Plane}. Recognizing
that each point on the grid represents a unique \tf\ shift
$\bs{H}_i$~(\ref{eq:TF_Basis_Element_H_i}) (with a corresponding
reflection coefficient~$s_i$), it is easy to see that every possible
target scene can be represented by some
matrix~$\bs{H}$~(\ref{eq:General_Unknown_Matrix_H}). If the number
of targets $K\ll N^2$, then the \tf\ grid will be sparsely
populated. By ``vectorizing'' the grid, we can represent it
as an {$N^2\!\times 1$} \textbf{sparse vector} $\bs{s}$. %Now the
%setup matches the general case of a sparse \rep\ under the \tf\
%basis introduced earlier.

Assume that the Alltop sequence $\bs{f}_\Alltop$ is sent by the
transmitter\footnote{The transmitter in Fig.~\ref{fig:Simple_Radar}
sends \emph{analog} signals. We assume here that there exists a
continuous signal which when discretized is the Alltop sequence
(\ref{eq:Alltop_sequnce}).}. The received signal now is of the form
in~(\ref{eq:Reformulation_Observation}). If the number of targets
obey the sparsity constraints in
Theorems~{\ref{thm:Matrix_Recovery_Guaranteed}-\ref{thm:Matrix_Recovery_Guaranteed_Noise}}
then we will be able to reconstruct the original target scene using
\cs\ techniques. Moreover, the resolution of the recovered target
scene is limited by how the \tf\ plane is discretized as dictated by
the~$N^2$ unique \tf\ shifts. That is, \emph{multiple targets
located at adjacent grid points can be resolved} due to the nature
of \cs\ reconstruction. The effect of discretization on the
resolution is discussed in more detail in the next section.

In reality, we are not actually ``beating'' the classical \up\ as
claimed above. Rather, we are just transferring to a different
mathematical perspective. The new \emph{\cs\ \up} is dictated by the
sparsity constraints of
Theorems~{\ref{thm:Matrix_Recovery_Guaranteed}-\ref{thm:Matrix_Recovery_Guaranteed_Noise}}.

It is interesting to note that Alltop specifically mentions the
applicability of his sequence to spread-spectrum radar. The cubic
phase in (\ref{eq:Alltop_sequnce}) is known in classical radar as a
discrete quadratic chirp, which is similar to what bats use to
``image'' their environment (although bats use a continuous
\emph{sonar} chirp). The use of a chirp is an effective way to
transmit a wide-bandwidth signal over a relatively short time
duration. However, here in \cs\ radar we make use of the incoherence
property of the Alltop sequence, which is due to specific properties
of prime numbers. Recall the three key points of this novel
approach: (1) the transmitted signal must be \textbf{incoherent},
(2) there is \textbf{no matched filter}, (3) instead, \cs\
techniques are used to recover the \textbf{sparse} target scene.

%%%%%%%%%%%%%%%%%%%%%%%%%%%%%%%%%%%%%%%%%%%%%%%%%%%%%%%%%%%%%%%%%%%%%%%%%%%%
\subsection{Comparison of Resolution Limits} %Resolution limit for CS Radar
%%%%%%%%%%%%%%%%%%%%%%%%%%%%%%%%%%%%%%%%%%%%%%%%%%%%%%%%%%%%%%%%%%%%%%%%%%%%

In this section we analyze the resolution limit for \cs\ radar and
compare it to the resolution limit dictated by the radar \up.

Assume that the transmitted signal is bandlimited to $[-\half B_1,
\half B_1]$. Actually, the received signal will have a somewhat
larger bandwidth $B > B_1$ due to the Doppler effect. However, in
practice this increase in bandwidth is small, so we can assume
$B\approx B_1$. We observe the signal over a duration\footnote{We
assume a periodic model here which can be relaxed using standard
zero-padding procedures.} $T$ and for simplicity sample it at the
Nyquist rate~$B$. That means we gather $N = BT$ many samples during
the observation interval. It is well-known that observing a signal
over a duration period~$T$ gives rise to a maximum frequency
resolution of $1/T$. The time resolution is equal to the Nyquist
sampling rate, i.e., $1/B$. The step-size for the discretization of
the \tf\ plane is therefore limited to $1/T$ and $1/B$,
respectively.

If $N\ge5$ is prime then we can use the Alltop sequence as described
in the previous section and recover multiple targets with a
resolution of $1/N$. Otherwise, there exist other ``incoherent''
sequences which can provide similar results to
Theorems~{\ref{thm:Matrix_Recovery_Guaranteed}-\ref{thm:Matrix_Recovery_Guaranteed_Noise}};
and therefore, can also achieve a resolution of $1/N$. \emph{Thus
for fixed~$T$ and fixed $B$, the smallest rectangle in the
time-frequency plane which can be resolved with \cs\ radar has size
$1/T \times 1/B = 1/N$.}\footnote{Note that a precise analysis on
the resolution limits of \cs\ radar must also take into account
approximating the continuous-time, continuous-frequency,
infinite-dimensional radar model by a discrete, finite-dimensional
model. We will report on this topic in a forthcoming paper.}

Now consider the Heisenberg uncertainty box associated with the
radar \up. When $\eps = 0$ in
(\ref{eq:Radar_Uncertainty_Princ_Conclusion}) this box must have an
\emph{area of at least unity}. This lower bound determines the resolution
limit of classical radar. Juxtapose this with the resolution limit
of \cs: we can easily make this box smaller by increasing the
observation period $T$ and/or the bandwidth $B$.\footnote{There
are, of course, practical considerations that prevent implementing
an extremely large observation period and/or bandwidth, which we
ignore for the purpose of this paper.} Therefore, in theory, \cs\ radar can
achieve better resolution than conventional radar.

Can we achieve an even better resolution than $1/N$ for fixed
duration $T$ and fixed bandwidth $B$ with compressed sensing radar?
Not with the existing theory and the existing algorithms. To achieve
better resolution one might be tempted to increase the sampling
rate. However, oversampling introduces correlations between the
samples, therefore it would not improve the incoherence of the
columns of $\bs{\Phi}$ (in practice though we always oversample
signals, but for different reasons).

The lower limit of $1/\gr{TB}$ appears in other areas of classical
radar as well, usually in the context of ``thumbtack'' functions. A
function is ``thumbtack-like'' if all of its values are close to
zero except for a unique large spike. These waveforms are also
sometimes referred to as ``low-correlation'' sequences. Due to
Properties~1 and 2 of the Alltop sequence in
Section~\ref{property:One} we see that its ambiguity surface
actually has this thumbtack feature too. Other thumbtack-like
ambiguity surfaces include those associated with the waveforms which
generate the equiangular line sets found in \cite{Calderbank}. The
crucial difference here is that, in general, the lower resolution
limit of $1/\gr{TB}$ can only be achieved in classical radar if
there is \emph{just one target}. As soon as several targets are
clustered together then interference from the non-zero portions of
the ambiguity function causes false positives. This dictates the
resolution limit, i.e., how close targets can be and still be able
to reliably distinguish them. The next section show computer
simulations which demonstrate this.

%%%%%%%%%%%%%%%%%%%%%%%%%%%%%%%%%%%%%%%%%%%%%%%%%%%%%%%%%%%%%%%%%%%%%%%%%%%%
\subsection{\CS\ and Classical Radar Simulations}
\begin{figure}[!ht]
\centering \subfigure{    \label{fig:subfig1:a}
    \includegraphics[width=1.6725in]{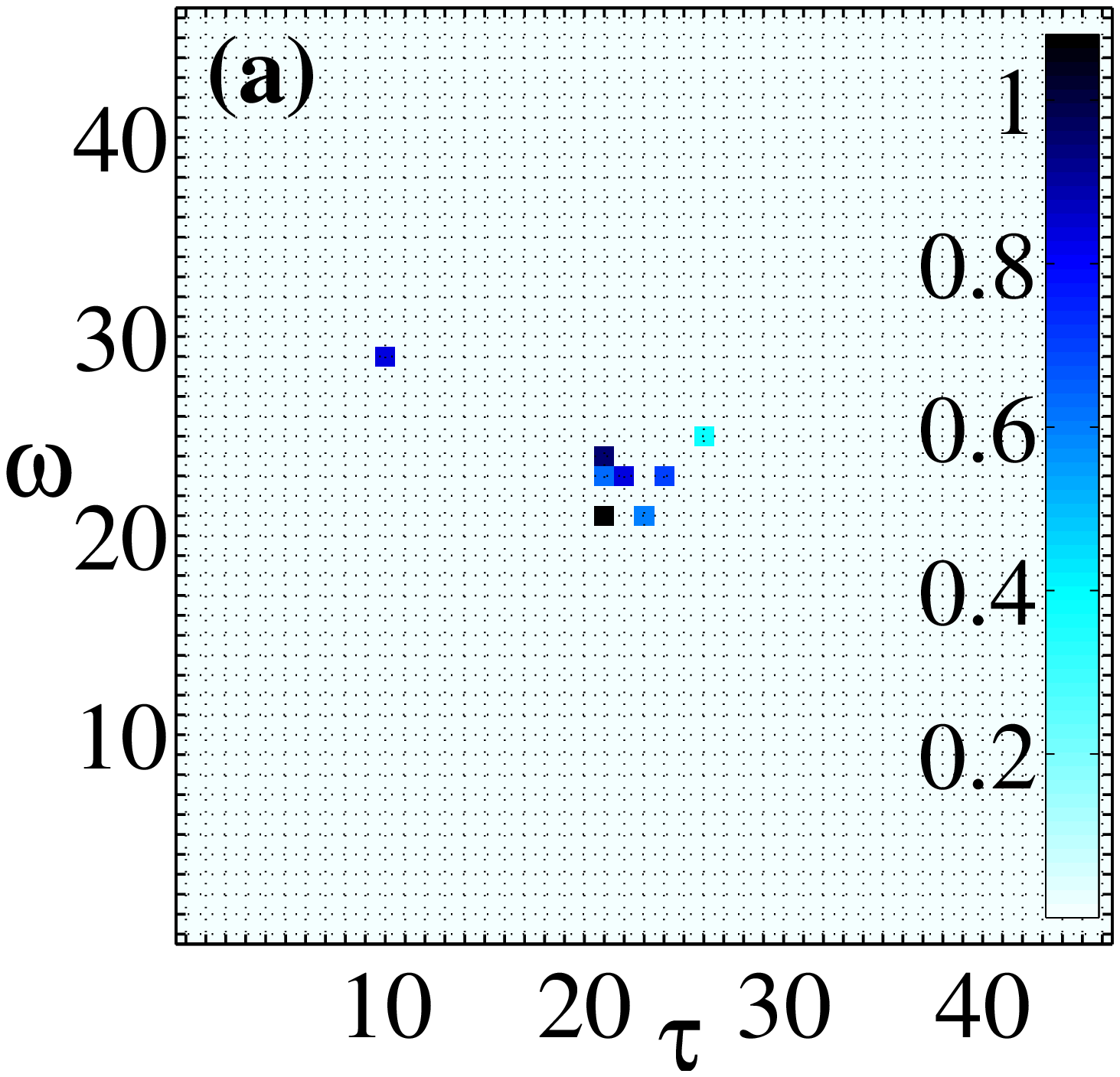}}
\subfigure{    \label{fig:subfig1:b}
    \includegraphics[width=1.6725in]{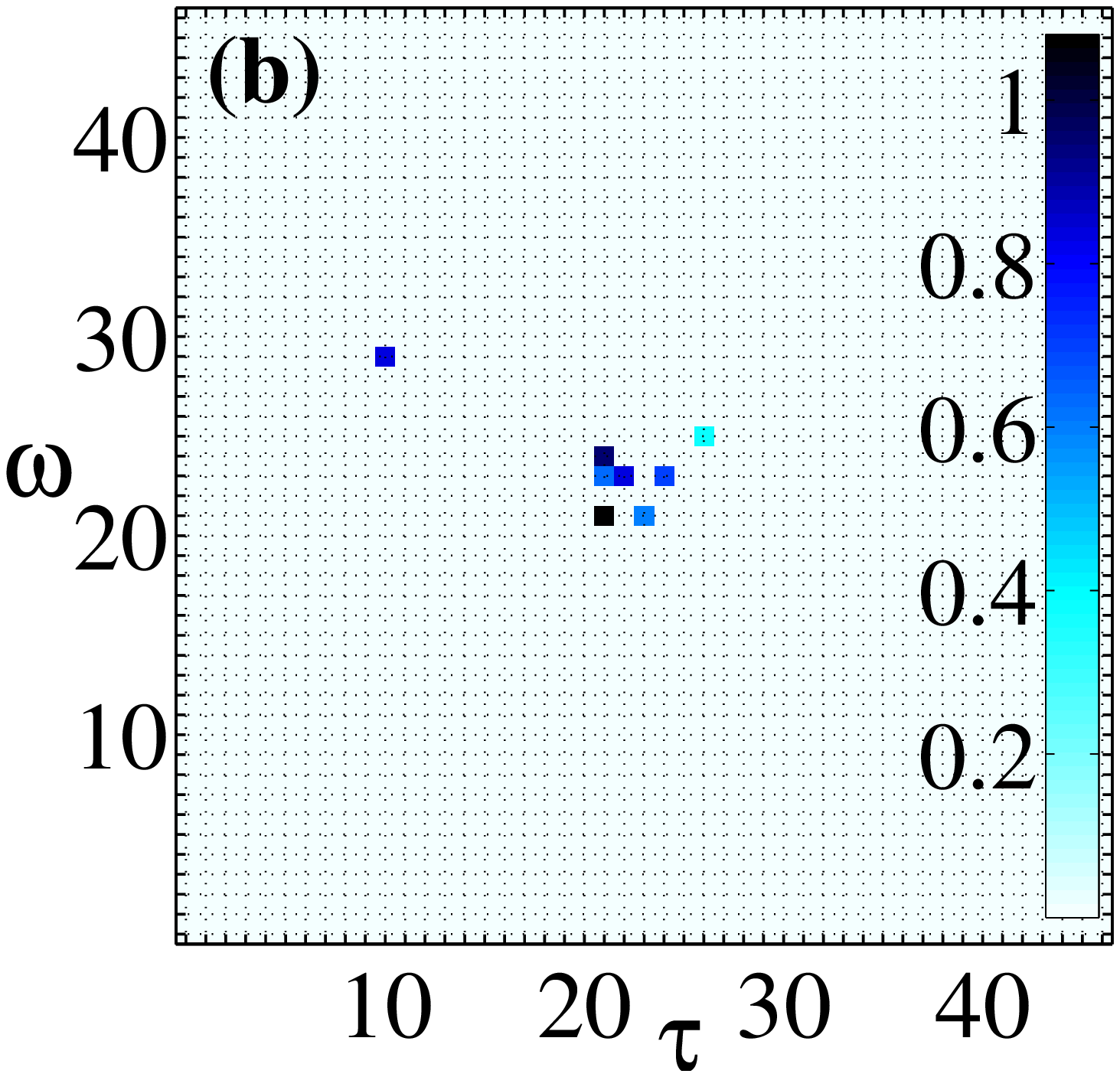}}\\
\subfigure{    \label{fig:subfig1:c}
    \includegraphics[width=1.6725in]{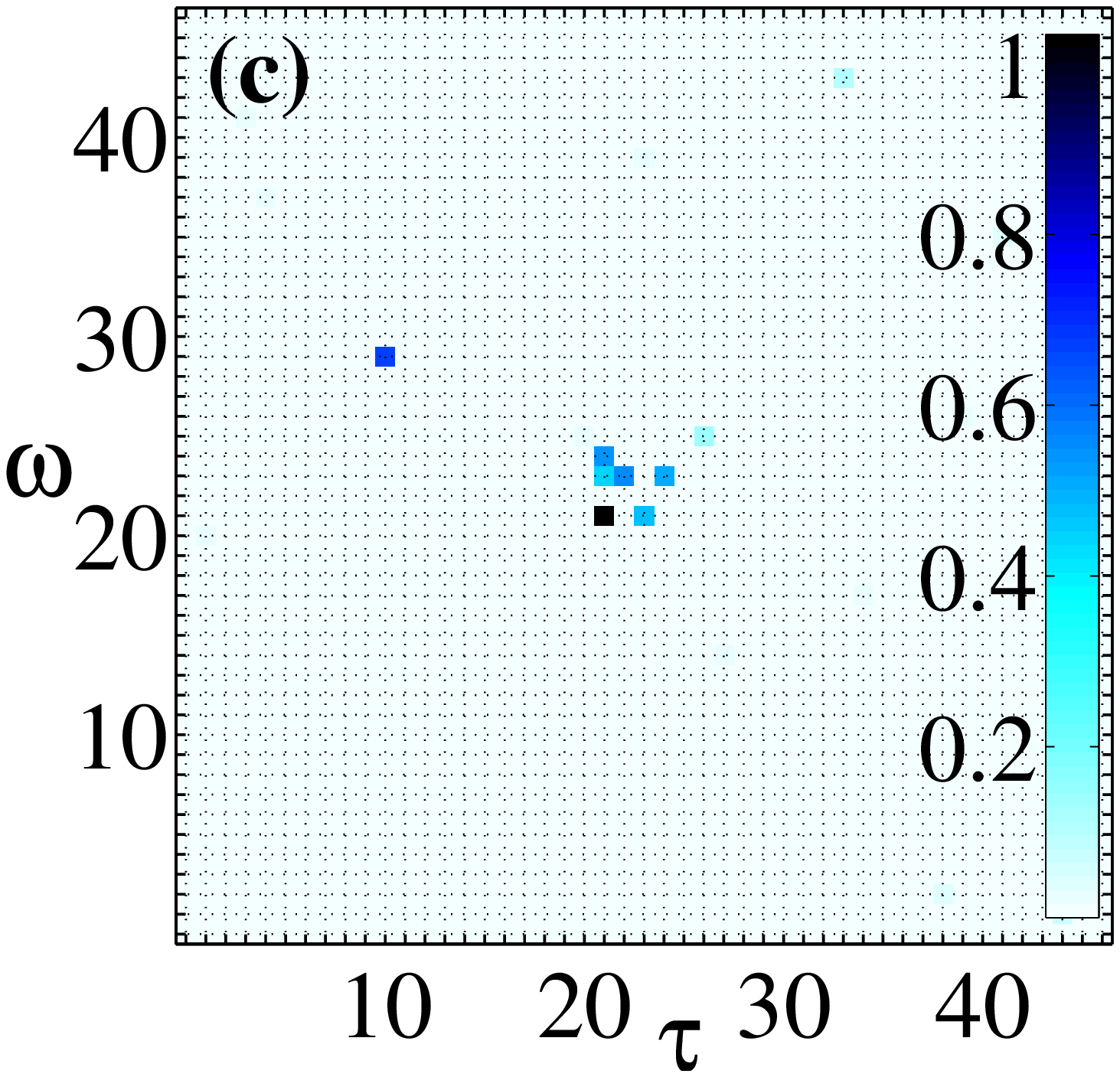}}
\subfigure{    \label{fig:subfig1:d}
    \includegraphics[width=1.6725in]{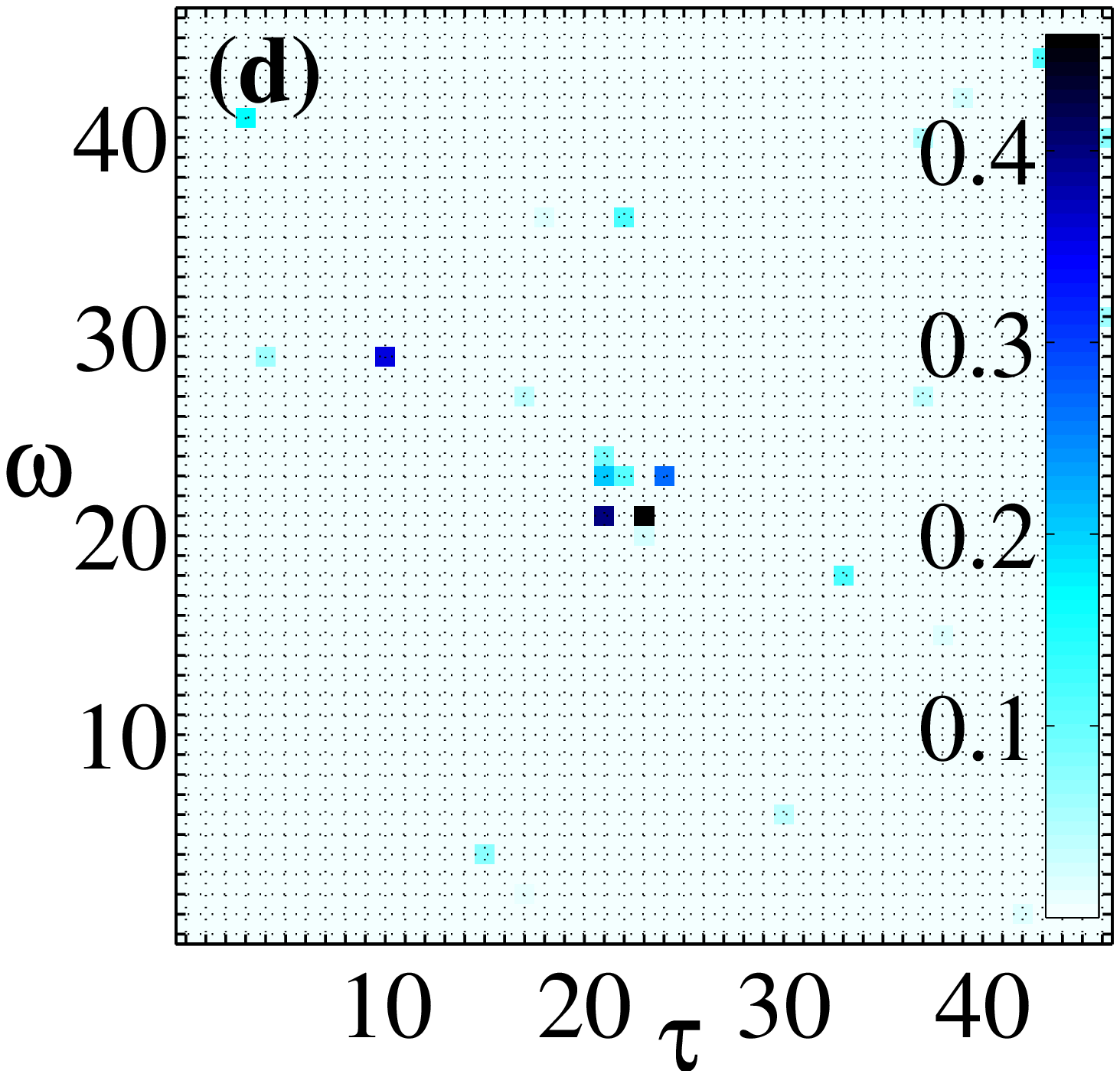}}
\caption{Radar simulation with $K=8$ targets on a $47\times47$ \tf\
grid. (a) Original target scene. \Cs\ reconstruction of original
target scene with SNR: (b) $\infty$ dB, (c) 15 dB, (d) 5 dB. Notice
in (b) that \cs\ perfectly recovers (a) in the case of no noise.}
\label{fig:First_Radar_Demo}
\end{figure}

\begin{figure}[!ht]
\centering \subfigure{    \label{fig:subfig2:a}
    \includegraphics[width=1.6725in]{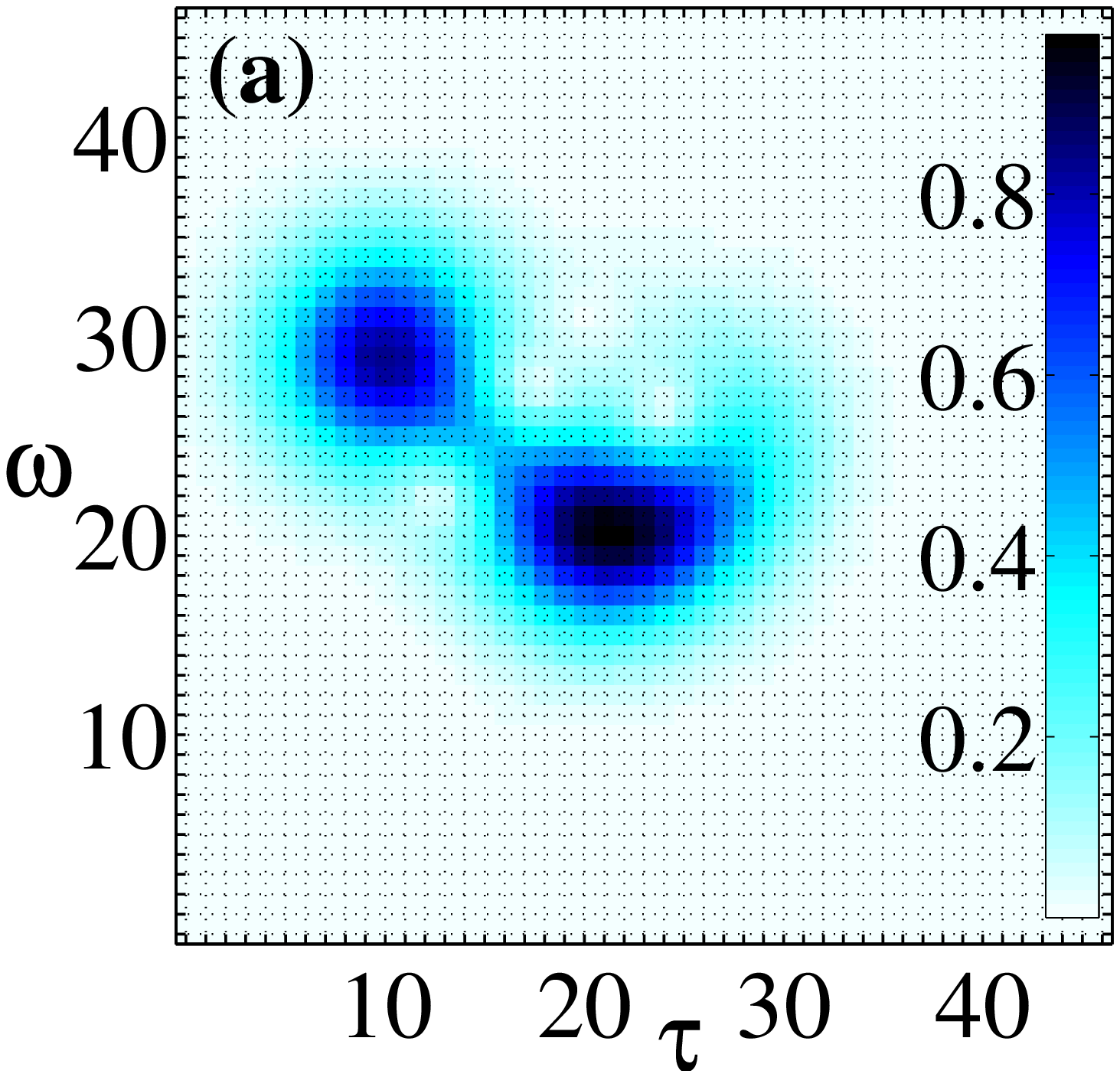}}
\subfigure{    \label{fig:subfig2:b}
    \includegraphics[width=1.6725in]{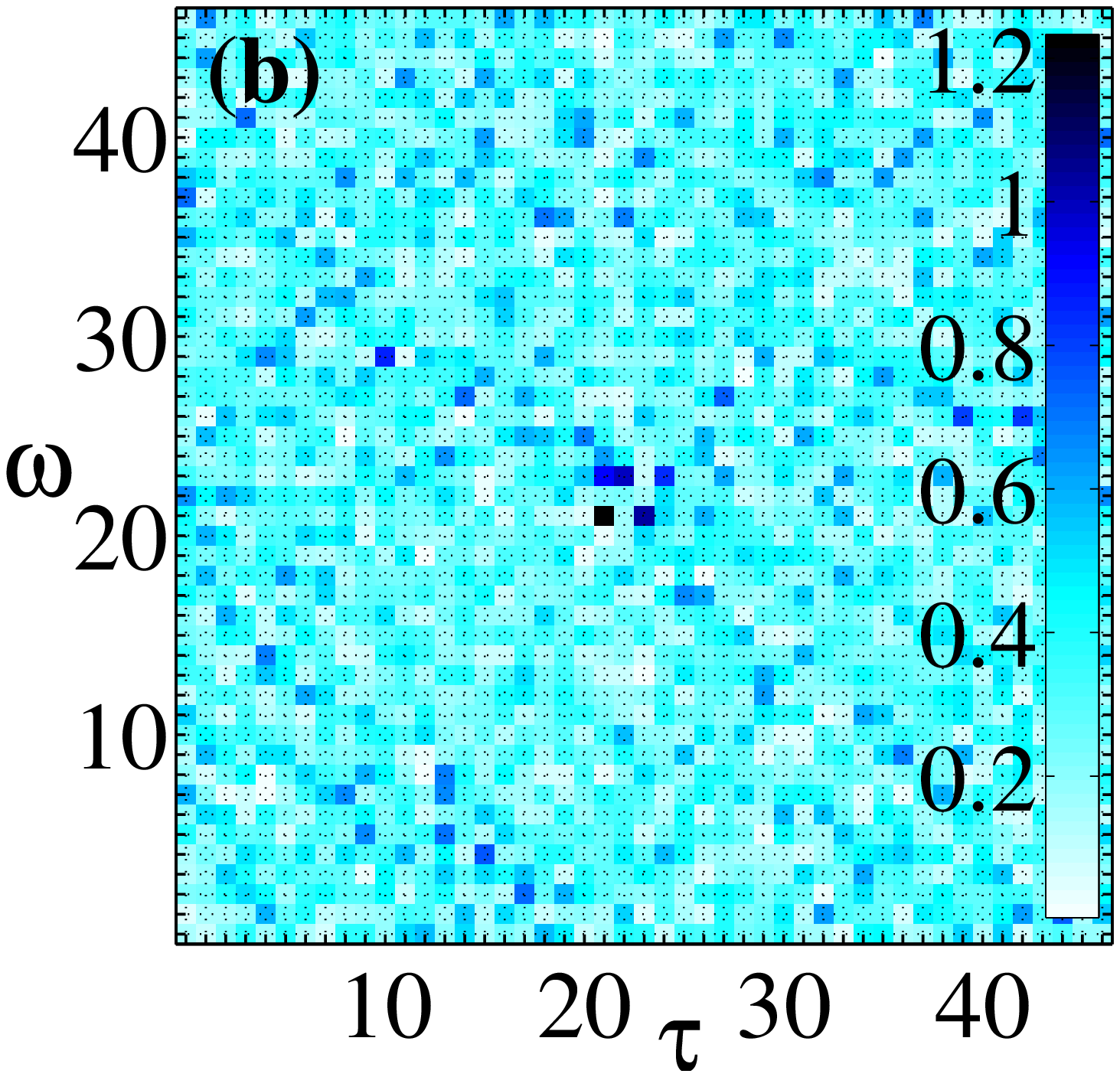}}\\
\subfigure{    \label{fig:subfig2:c}
    \includegraphics[width=1.6725in]{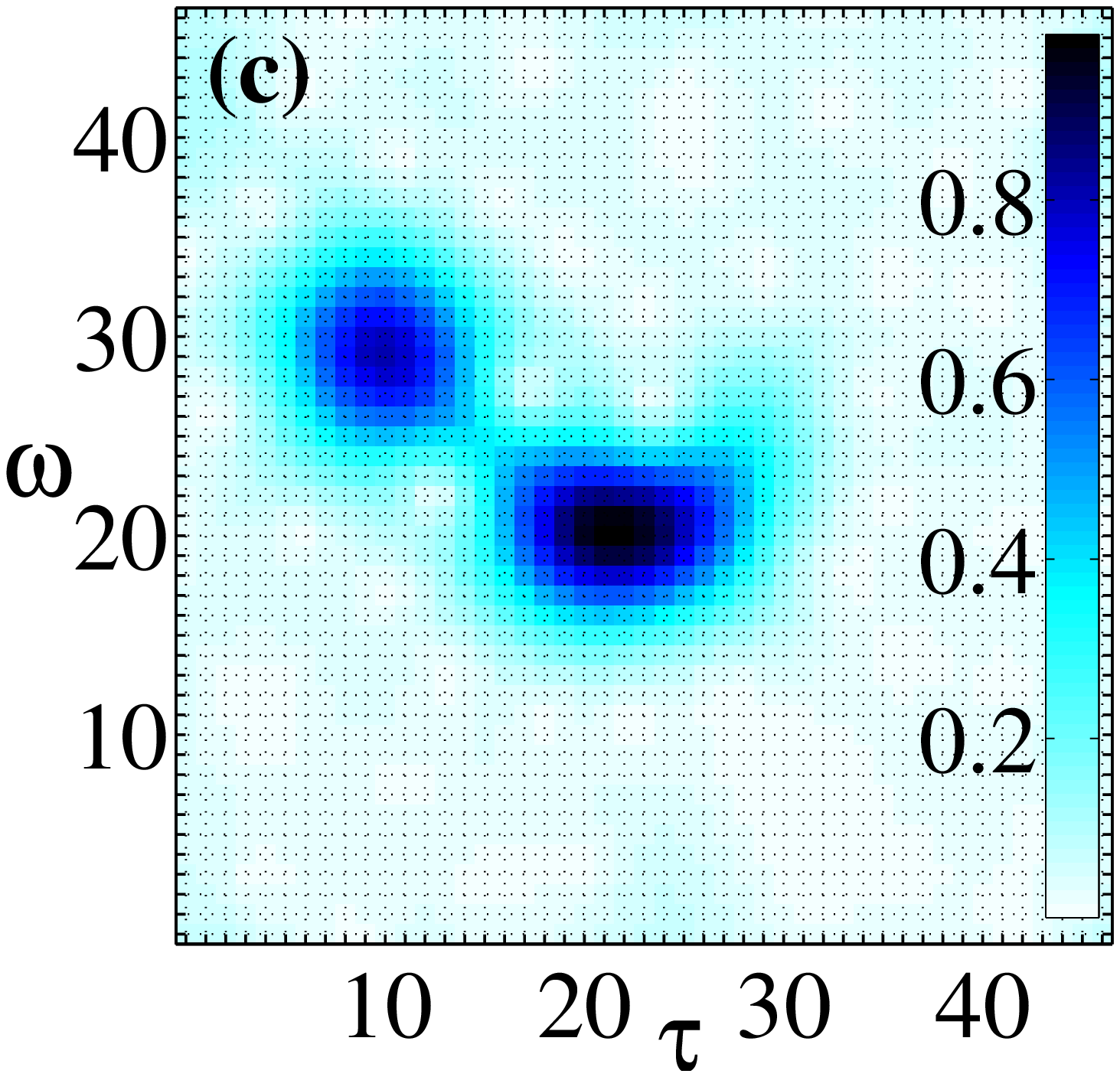}}
\subfigure{    \label{fig:subfig2:d}
    \includegraphics[width=1.6725in]{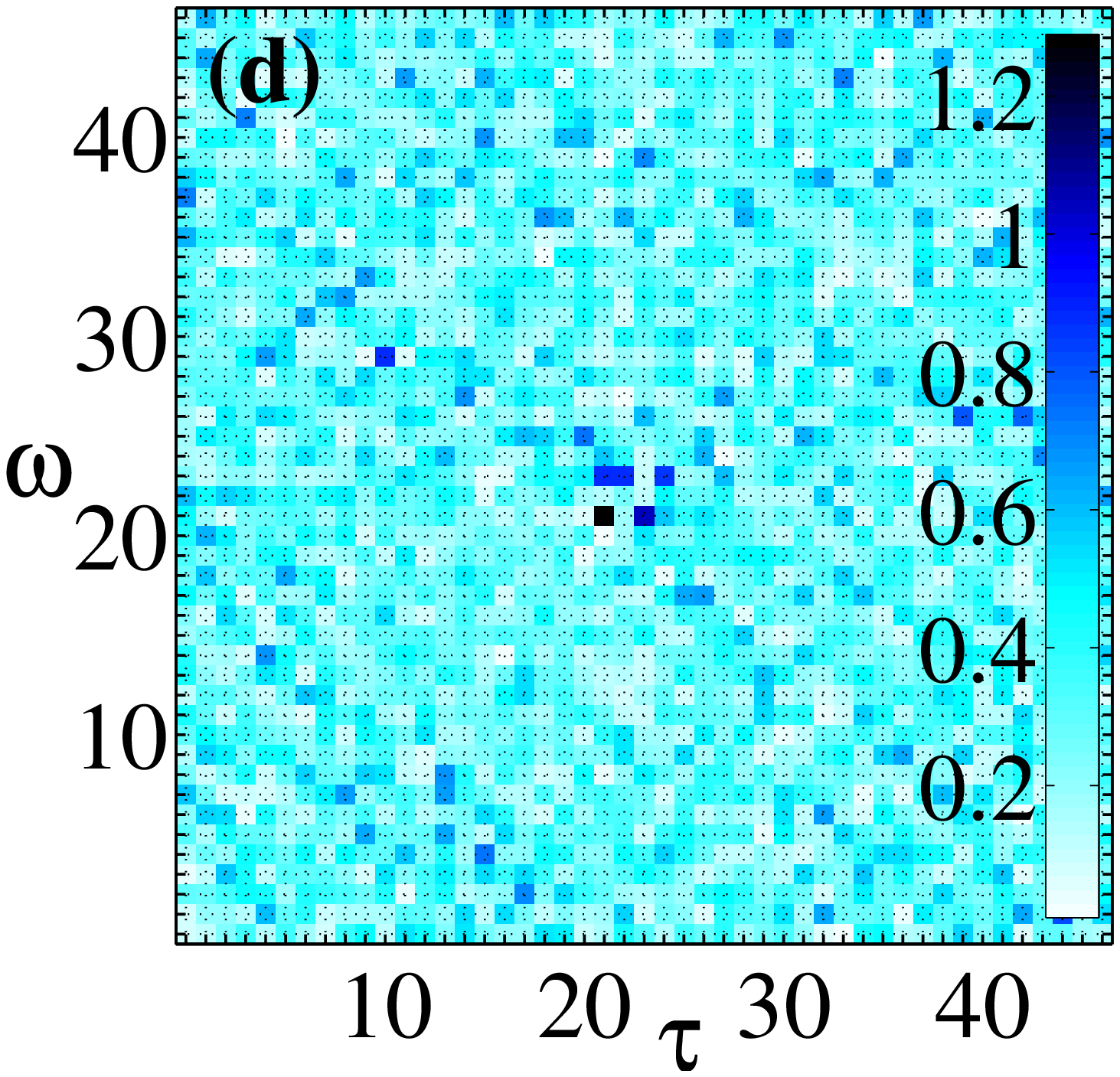}}\\
\subfigure{    \label{fig:subfig2:e}
    \includegraphics[width=1.6725in]{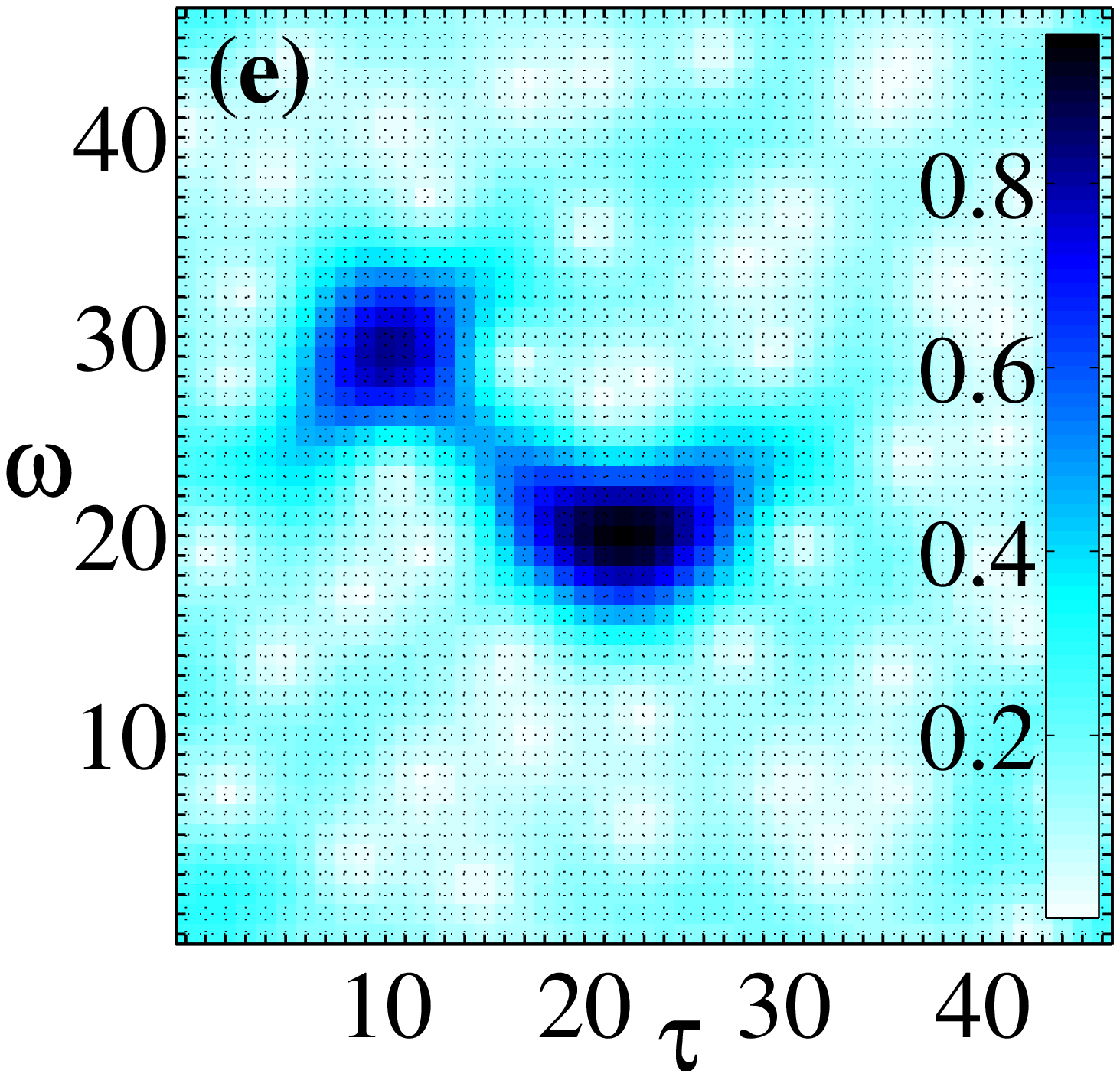}}
\subfigure{    \label{fig:subfig2:f}
    \includegraphics[width=1.6725in]{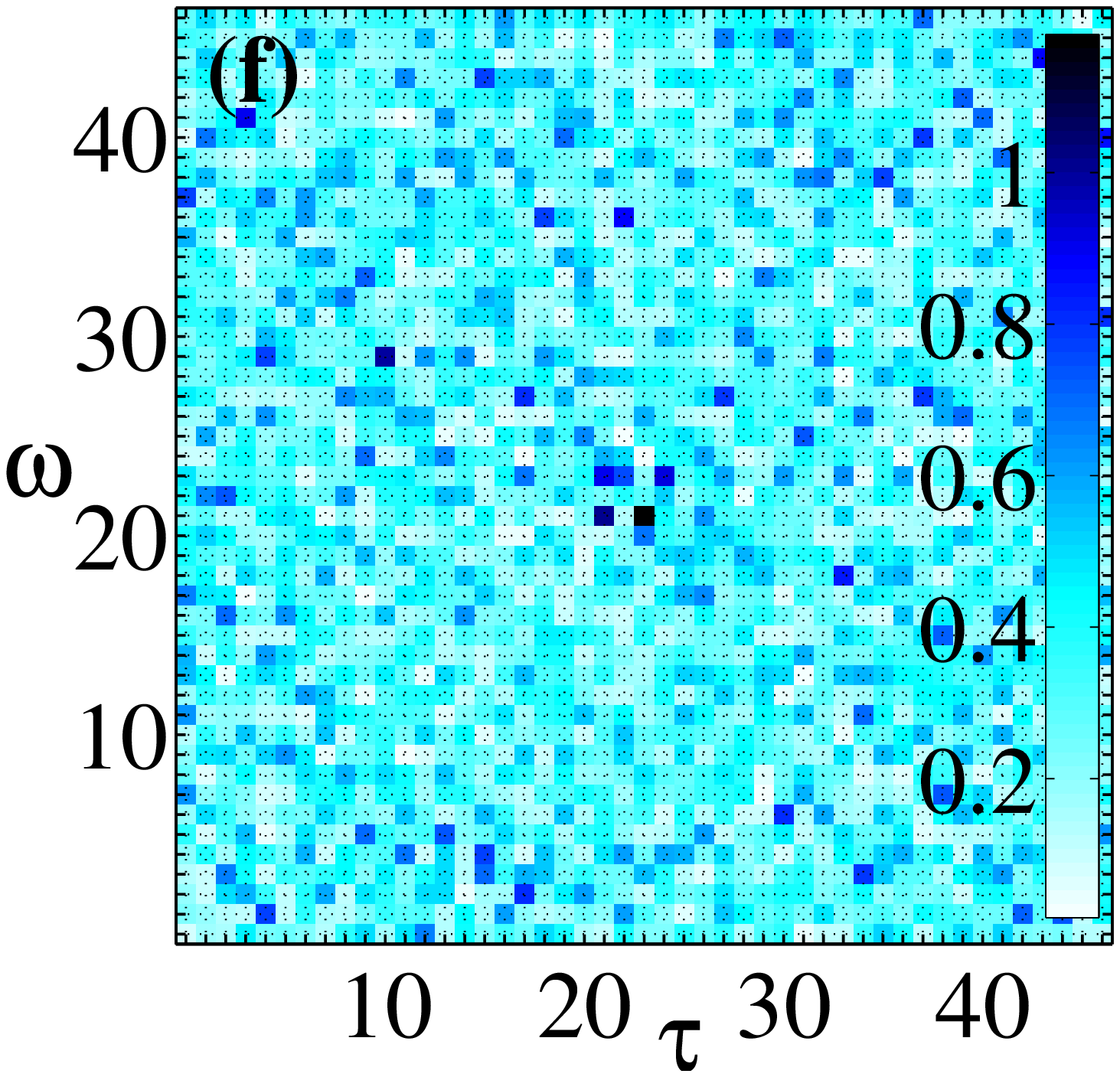}}
\caption{Traditional radar reconstruction of
Fig.~\ref{fig:subfig1:a}'s original target scene. With no noise: (a)
Gaussian pulse, (b)~Alltop sequence. With SNR = 15 dB: (c) Gaussian
pulse, (d) Alltop sequence. With SNR = 5 dB: (e) Gaussian pulse, (f)
Alltop sequence.}
\label{fig:Second_Radar_Demo}
\end{figure}

%%%%%%%%%%%%%%%%%%%%%%%%%%%%%%%%%%%%%%%%%%%%%%%%%%%%%%%%%%%%%%%%%%%%%%%%%%%%
Figures~\ref{fig:First_Radar_Demo} and \ref{fig:Second_Radar_Demo}
show the result of Matlab radar simulations. For purposes of
normalization the grid spacing in these figures is ${1}/{\sqrt{N}}$.
Hence, the numbers shown on the axes represent multiples of
${1}/{\sqrt{N}}$. A random \tf\ scene with $K=8$ targets and $N=47$
is presented in Figure~\ref{fig:subfig1:a}. The \emph{\cs\ radar}
simulation used the Alltop sequence to identify the targets. In the
noise-free case of Figure~\ref{fig:subfig1:b} it is clear that \cs\
was able to \emph{perfectly reconstruct} the target scene
($\normtwo{\bs{s}-\bs{s^\bigstar}} \sim 10^{-8}$). Moreover, it is
obvious that targets located at \emph{adjacent grid points} can be
resolved, confirming the discussion of the last section.

Figure~\ref{fig:subfig1:c} shows how \cs\ starts to suffer in the
presence of additive white Gaussian noise (AWGN). Here the
signal-to-noise ratio (SNR) is 15 dB. Some faint false positives
have appeared, yet the target scene has
still been identified. The performance %of \cs\ radar
with 5 dB SNR is shown in Figure~\ref{fig:subfig1:d}. One target was
lost, many false positives have appeared, and the magnitudes of the
targets have been significantly reduced. Clearly, these are all
undesirable effects. It remains an open problem in the \cs\
community how to deal with such noisy situations.

As a comparison to \cs\ Figure~\ref{fig:Second_Radar_Demo} presents
\emph{classical radar} reconstruction (which \textbf{uses a matched filter}
as described in Section~\ref{sect:Classical_Radar}) with two different
transmitted pulses. The ambiguity surfaces associated with these two
waveforms demonstrate, in some sense, two extremes of traditional radar
performance. In the first case, the ambiguity surface is a
relatively wide Gaussian pulse, whereas in the second case the
ambiguity surface is a highly concentrated ``thumbtack'' function.
We stress that these are not necessarily the final results of
traditional target reconstruction, and are included only for rough
comparison. In practice, radar engineers use extremely advanced
techniques to determine target range and velocity.

Figures~\ref{fig:subfig2:a}, \ref{fig:subfig2:c}, and
\ref{fig:subfig2:e} show the original target scene of
Figure~\ref{fig:subfig1:a} reconstructed using a Gaussian pulse. The
(self) ambiguity function associated with a Gaussian pulse is a
two-dimensional (2D) Gaussian pulse as a result of the STFT
in~(\ref{eq:Cross_Ambiguity_equiv_STFT}). Therefore, according to
(\ref{eq:Matched_Filter_Ambiguity_Eqn}) we see that the radar scenes
in these figures consist of a 2D Gaussian pulse centered at each
target in the \tf\ plane. In each of these it is clear that the
targets in the center are contained within the Heisenberg boxes of
its neighbors. Depending on the sophistication of subsequent
algorithms some of the targets may be unresolvable. It is also
apparent that Figures~\ref{fig:subfig2:c} and~\ref{fig:subfig2:e}
suffer from added noise, and this compounds the problem of accurate
resolution \cite{Rihk}.

As a consequence of the grid spacing, the Heisenberg box associated
with the Gaussian pulse's ambiguity surface has been normalized to a
square of unit area. This roughly corresponds to the support size of
$U$ in (\ref{eq:Radar_Uncertainty_Princ_Conclusion}), and is
empirically verified in Figure~\ref{fig:subfig2:a} where we see that
the diameter of the uncertainty region around the isolated target at
$\gr{\tau,\omega}=\gr{10,29}$ spans approximately seven grid points.
Since the grid spacing is~${1}/{\sqrt{N}}$ we confirm that the base
and height of the Heisenberg box are each approximately
${7}/{\sqrt{47}}\approx 1$.

Returning to the discussion of the previous section, it is clear
that the noise-free cases shown in Figures~\ref{fig:subfig1:b} and
\ref{fig:subfig2:a} experimentally confirm that \cs\ radar can
achieve much higher resolution than traditional
techniques.\footnote{ There are many different ways to determine
resolution in classical radar. Moreover, in the presence of noise,
the SNR must also be incorporated. See~\cite{Rihk,IMA}.} To make the
comparison fair, we are using the same number of observations in the
recovery for both \cs\ and classical radar. In this sense, it
becomes apparent that we are leveraging the power of \cs\ theory in
a different way than explained in Section~\ref{sect:CS}. The typical
\cs\ application makes \emph{far fewer} observations than necessary
and still obtains perfect reconstruction of the data. However, in
this model of \cs\ radar we implicitly assume Nyquist sampling of
the baseband signal. Therefore, with this setup, the benefit of
employing \cs\ recovery manifests itself as a dramatic
\emph{increase in resolution}.

In contrast with a Gaussian pulse we now examine a waveform whose
associated ambiguity surface is thumbtack-like.
Figures~\ref{fig:subfig2:b}, \ref{fig:subfig2:d}, and
\ref{fig:subfig2:f} depict the original target scene
\emph{traditionally reconstructed} using the Alltop sequence. Take
note of the distinction with \cs\ radar presented in
Section~\ref{sect:CS_Radar} which also uses this function. Here, the
classical approach transmits the Alltop sequence, and then uses a
\textbf{matched filter} to correlate the received signal with a \tf\
shifted Alltop sequence as
in~(\ref{eq:Matched_Filter_Ambiguity_Eqn}). The radar scene will now
consist of a thumbtack function centered at each target. In theory,
this radar would provide target resolution similar to our \cs\
version (i.e., the target is represented as a point source in \tf\
plane rather than a ``spread out'' uncertainty region).

However, the situation is not so simple. The non-zero portions of
the ambiguity function can accumulate to create undesirable effects.
This is shown in Figure~\ref{fig:subfig2:b} where it is apparent,
\emph{even in the ideal case of no added noise}, that there is a
great deal of interference. Moreover, this type of ``noise'' is
deterministic and cannot be remedied by averaging over multiple
observations. Notice that the interference seems to be distributed
over a wide range of amplitudes. In fact, referring to the original
target scene in Figure~\ref{fig:subfig1:a}, it appears that some of
the weaker targets (i.e., the ones with the smallest reflection
coefficient in magnitude) have been buried in this noise. Even if a
reasonable threshold could be determined, perhaps only a few of the
strongest targets would be detected and many false positives would
remain. This is a substantial problem since the dynamic range of the
targets can be quite large.

We present these results to emphasize that naive application of
\emph{traditional radar} techniques with the Alltop sequence will
fail if the radar scene contains more than just a few strong
targets. The outcome will be similar if other low-correlation
sequences are used.

Regardless of whether a transmitted waveform has an ambiguity
surface which is spread or narrow, interference from adjacent
targets will necessarily occur in classical radar, and this will
result in undesirable effects. In contrast, \cs\ radar does not
experience this interference since it completely dispenses with the
need for a matched filter. Therefore, there are no issues with the
ambiguity function of the transmitted signal.

%%%%%%%%%%%%%%%%%%%%%%%%%%%%%%%%%%%%%%%%%%%%%%%%%%%%%%%%%%%%%%%%%%%%%%%%%%%%
%%%%%%%%%%%%%%%%%%%%%%%%%%%%%%%%%%%%%%%%%%%%%%%%%%%%%%%%%%%%%%%%%%%%%%%%%%%%
\section{Other Applications} \label{sect:Other_Applications}
%%%%%%%%%%%%%%%%%%%%%%%%%%%%%%%%%%%%%%%%%%%%%%%%%%%%%%%%%%%%%%%%%%%%%%%%%%%%
%%%%%%%%%%%%%%%%%%%%%%%%%%%%%%%%%%%%%%%%%%%%%%%%%%%%%%%%%%%%%%%%%%%%%%%%%%%%
\emph{Narrowband} radar is by no means the only application to which
the techniques presented here can be used. \emph{Wideband} radar
systems admit a received signal which is of the form
$$\abs{r(t)} \;\sim\; \Abs{f\Gr{\frac{t-a\:\!\tau_x}{a}}},
        \qquad a = 1-2v/c.$$
This shift-scaled signal is well-represented by a wavelet basis, and
it seems feasible to replace the time-frequency dictionary by a
properly chosen time-scale dictionary. In a different direction, the
methods introduced in this paper can also be extended to
multiple-input multiple-output (MIMO) radar systems.

Our approach can also be applied, with suitable modifications, to
other applications that involve the identification of a linear
(time-varying) system. For instance, a challenging task in
underwater acoustic communication is the estimation of the acoustic
propagation channel. Unlike mobile radio channels, underwater
acoustic channels often exhibit large delay spreads with substantial
Doppler shifts. Of course, the location of the scatterers and the
amount of Doppler shift are \emph{a priori} not known. However, it is known
that underwater communication channels do have a sparse
representation in the time-frequency domain, e.g.,
see~\cite{LiPreisig}. Thus, there is a good chance that our approach
via compressed sensing can lead to a channel estimation method that
provides higher resolution than conventional methods. We point out
that in order to turn compressed sensing-based underwater acoustic
channel estimation into a reliable method, one needs to carefully
incorporate various other properties of underwater environments,
e.g., whether we are dealing with a deep sea environment or a
shallow water environment.

Another application where the proposed compressed sensing approach
seems useful arises in high-resolution radar imaging. For instance,
when we consider the imaging of (moving) point targets, one would
need to combine our time-frequency based approach with the Born
approximation of Helmholtz's equation. This approach is a topic of
our current research.

Other applications arise in blind source
separation~\cite{ShanSwaryAviyente}, sonar, as well as underwater
acoustic imaging based on matched field processing.

%%%%%%%%%%%%%%%%%%%%%%%%%%%%%%%%%%%%%%%%%%%%%%%%%%%%%%%%%%%%%%%%%%%%%%%%%%%%
%%%%%%%%%%%%%%%%%%%%%%%%%%%%%%%%%%%%%%%%%%%%%%%%%%%%%%%%%%%%%%%%%%%%%%%%%%%%
\section{Discussion}
%%%%%%%%%%%%%%%%%%%%%%%%%%%%%%%%%%%%%%%%%%%%%%%%%%%%%%%%%%%%%%%%%%%%%%%%%%%%
%%%%%%%%%%%%%%%%%%%%%%%%%%%%%%%%%%%%%%%%%%%%%%%%%%%%%%%%%%%%%%%%%%%%%%%%%%%%
We have provided a sketch for a high-resolution radar system based
on \cs. Assuming that the number of targets obey the sparsity
constraint in Theorem~\ref{thm:Matrix_Recovery_High_Prob}, the
Alltop sequence can perfectly identify the radar scene with high
probability using \cs\ techniques. Numerical simulations confirm
that this sparsity constraint is too strict and can be relaxed to $K
\leq N/\gr{2\log{N}}$, although this has yet to be proven mathematically.

It must be emphasized that our model presents radar in a rather
simplified manner. In reality, radar engineers employ highly
sophisticated methods to identify targets. For example, rather than
a single pulse, a signal with multiple pulses is often used and
information is averaged over several observations. We also did not
address how to discretize the analog signals used in both \cs\ and
classical radar. A more detailed study covering these issues is the
topic of another paper.

Related to the discretization issue is the fact that \cs\ radar does
not use a matched filter at the receiver. This will directly impact
A/D conversion, and has the potential to reduce the overall data
rate and to simplify hardware design. These matters are discussed in
\cite{Baraniuk}, although it does not consider the case of moving
targets. In our study the major benefit of relinquishing the matched
filter is to avoid the target uncertainty and interference resulting
from the ambiguity function.

Since many of the implementation details of our \cs\ radar have yet
to be determined, and since classical radar can also be implemented
in many ways we were only able to make a rough comparison between
their respective resolutions. \emph{Regardless, the radar \up\ lies
at the core of traditional approaches and limits their performance.}
We contend that \cs\ provides the potential to achieve higher
resolution between targets. The radar simulations presented confirm
this claim.

It must be stressed again that the success of this stylized \cs\
radar relied on the \textbf{incoherence} of the
dictionary~$\bs{\Phi}_\Alltop$ resulting from the Alltop sequence.
There exist other probing functions with similar incoherence
properties. Numerical simulations with $\bs{f}$ as a random Gaussian
signal, as well as a constant-envelope random-phase signal indicate
similar behavior to what we have reported for the Alltop sequence.
At the time of writing this paper we became aware of a similar
study~\cite{PfanRauhTan} where the properties of these functions are
analyzed in the context of abstract system identification using \cs.

\begin{figure}[!h]
\centering \subfigure{    \label{fig:CounterExample_a}
    \includegraphics[width=1.6725in]
    {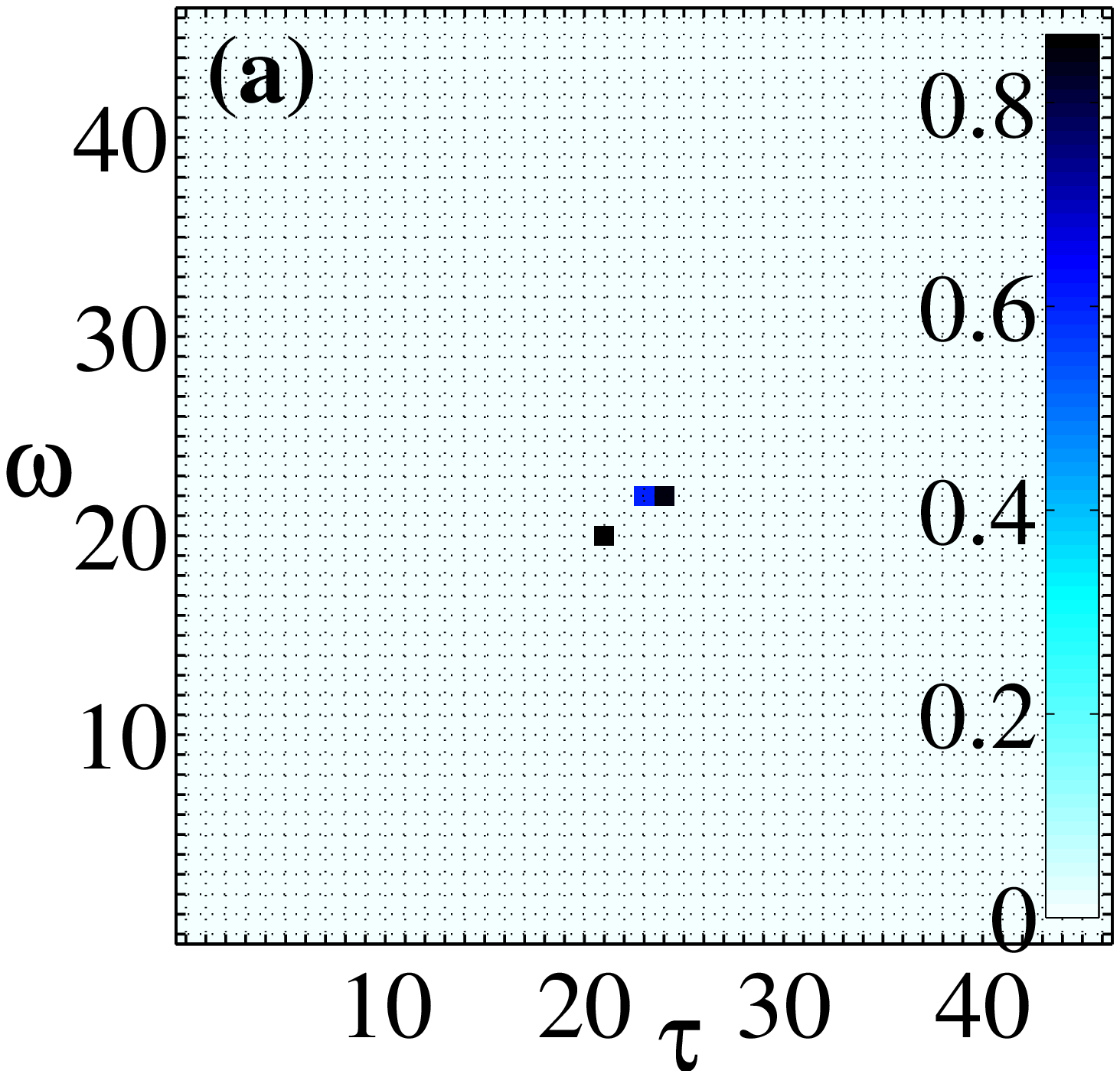}}
\subfigure{    \label{fig:CounterExample_b}
    \includegraphics[width=1.6725in]
    {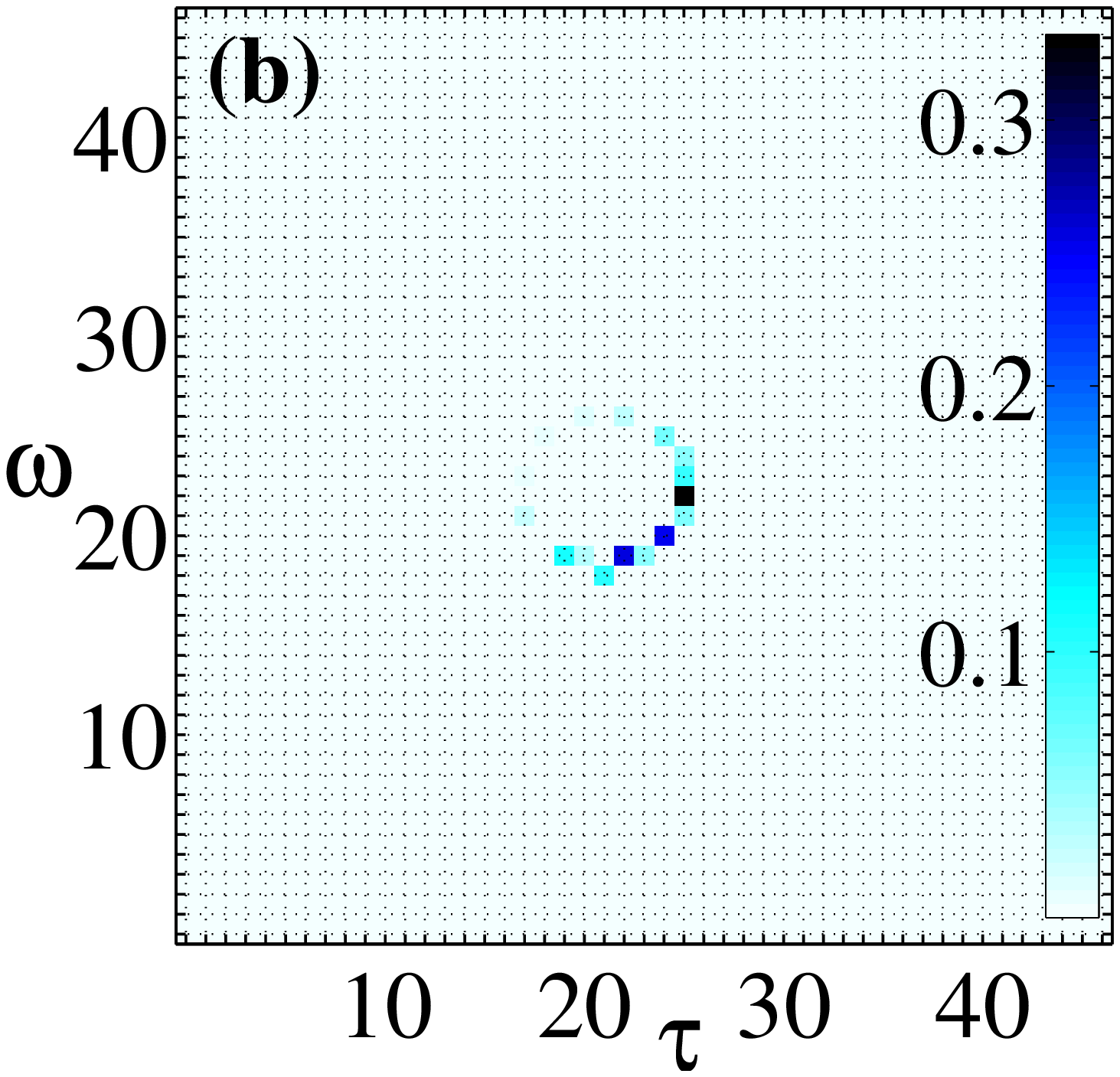}}
\caption{Radar simulation with $K=3$ targets on a $47\times47$ \tf\
grid. (a) Original target scene. (b) Traditional Gaussian pulse and
reconstruction using $\ell_1$ minimization (no noise). It is clear
that conventional radar with $\ell_1$ minimization completely fails.
However, Theorem~\ref{thm:Matrix_Recovery_Guaranteed}
\textbf{guarantees} perfect recovery in this case.}
\label{fig:Second_CounterExample}
\end{figure}

There is also the possibility of combining classical radar
techniques with $\ell_1$ recovery. Initial tests show that while we
get good reconstruction, the results are not guaranteed, even in the
case of no noise. Figure \ref{fig:Second_CounterExample} shows a
striking example. In this noise-free scenario, a Gaussian pulse
has been transmitted and reconstruction is done using $\ell_1$
minimization. Figure \ref{fig:CounterExample_a} shows an original
radar scene with $K=3$ targets. It is clear from Figure
\ref{fig:CounterExample_b} that \emph{none} of the targets
have been correctly recovered. In contrast,
Theorem~\ref{thm:Matrix_Recovery_Guaranteed} proves that we are
\emph{guaranteed} to perfectly recover both of these target scenes
when transmitting the Alltop sequence. (Note, in order to employ
Theorem~\ref{thm:Matrix_Recovery_Guaranteed}, we need to satisfy
$K<\half\gr{\sqrt{N}+1}$. With $N=47$ we can only use $K=3$ targets
since $3<\half\gr{\sqrt{47}+1}\approx3.93$.)

%%%%%%%%%%%%%%%%%%%%%%%%%%%%%%%%%%%%%%%%%%%%%%%%%%%%%%%%%%%%%%%%%%%%%%%%%%%%
%%%%%%%%%%%%%%%%%%%%%%%%%%%%%%%%%%%%%%%%%%%%%%%%%%%%%%%%%%%%%%%%%%%%%%%%%%%%
\appendices{}
\section{Proof of the Theorems \label{app:Proofs}}
%%%%%%%%%%%%%%%%%%%%%%%%%%%%%%%%%%%%%%%%%%%%%%%%%%%%%%%%%%%%%%%%%%%%%%%%%%%%
%%%%%%%%%%%%%%%%%%%%%%%%%%%%%%%%%%%%%%%%%%%%%%%%%%%%%%%%%%%%%%%%%%%%%%%%%%%%
For notational simplicity denote the coherence of dictionary
$\bs{\Phi}$ as $\mu$. We need the following theorems which deal with
incoherent dictionaries such as $\bs{\Phi}_\Alltop\in\C^{N\times
N^2}$. Recall for $\bs{\Phi}_\Alltop$ that $\mu =
{1}/{\sqrt{N}}$ with prime~$N\geq5$.

\begin{prop}[\cite{Tropp_Dictionaries}, Theorem~B] \label{prop:Tropp_Thm_B}
Let $\bs{X}$ be a random $K$-column subdictionary of $\bs{\Phi}$
(i.e, every $K$-column subset of $\bs{\Phi}$ has an equal
probability of being chosen). The condition
$\sqrt{\mu^2K\log{K}\cdot\vartheta} +
\frac{K}{N^2}\norm{\bs{\Phi}}^2 \leq c\delta$ with $\vartheta\geq1$
implies that $\Pr{\set{\norm{\bs{X^*X} - \bs{I}} \geq \delta}} \leq
K^{-\vartheta}$ where $c$ is an absolute constant.
\end{prop}

\begin{prop} [\cite{Tropp_Dictionaries}, Theorem~14] \label{prop:Tropp_Thm_14}
Suppose random $\bs{s}\in\C^{N^2}$ has support~$T$, sparseness $K =
\abs{T}$, and nonzero coefficients whose phases are uniformly
distributed on the interval $[0,2\pi)$. Set $\bs{y} \;\!=\;\!
\bs{\Phi}\bs{s}$, and let $\bs{\Phi}_T$ be the submatrix consisting
of the columns $\bs{\varphi}_j$ of $\bs{\Phi}$ for $j\in T$. Suppose
$8\mu^2K\leq1/\log{\gr{N^2\!/\zeta}}$ and that the least singular
value $\sigma_{\min}\gr{\bs{\Phi}_T}\geq 1/\sqrt{2}$. Then $\bs{s}$
is the unique solution to BP except with probability $2\zeta$.
\end{prop}

\begin{prop} [\cite{DonElad}, Theorem~3] \label{prop:DonElad_Thm_3}
Suppose a noisy signal $\bs{y} = \bs{\Phi}\bs{s} + \bs{e}$ is
constructed as a sparse combination of the columns of dictionary
$\bs{\Phi}\in\C^{N\times N^2}$ with coherence $\mu$. Assume the
sparsity of~$\bs{s}$ obeys $K < \gr{1+\mu}/\gr{2\mu + 4\eps
\sqrt{N}/T}$, and the entries of the noise are bounded
$\abs{e_n}\le\eps$. Then the solution~$\bs{s^\bigstar}$ to BP
exhibits stability $\normone{\bs{s}-\bs{s^\bigstar}}\le T$.
\end{prop}

%%%%%%%%%%%%%%%%%%%%%%%%%%%%%%%%%%%%%%%%%%%%%%%%%%%%%%%%%%%%%%%%%%%%%%%%%%%%
\subsection{Theorem~\ref{thm:Matrix_Recovery_Guaranteed}}
%%%%%%%%%%%%%%%%%%%%%%%%%%%%%%%%%%%%%%%%%%%%%%%%%%%%%%%%%%%%%%%%%%%%%%%%%%%%
\begin{IEEEproof}
Theorem B in \cite{Tropp_Greed} (which incorporates results from
\cite{DonElad}, \cite{elad02generalized}, and
\cite{GribonvalNielsen}) concludes for general
dictionary~$\bs{\Phi}$ that every $K$-sparse signal $\bs{s}$ with
$K< \half\gr{\mu^{-1}+1}$ is the unique sparsest representation, and
is guaranteed to be recovered by both BP and OMP when observing
$\bs{y} = \bs{\Phi}\bs{s}$. Set $\bs{\Phi} = \bs{\Phi}_\Alltop$ and
assume the hypothesis of
Theorem~\ref{thm:Matrix_Recovery_Guaranteed}.
Equation~(\ref{eq:Reformulation_Observation}) provides $\bs{y} =
\bs{H}\bs{f}_\Alltop = \bs{\Phi}_\Alltop\bs{s}$. The result follows
by substituting ${\mu = 1/{\sqrt{N}}}$.
\end{IEEEproof}

%%%%%%%%%%%%%%%%%%%%%%%%%%%%%%%%%%%%%%%%%%%%%%%%%%%%%%%%%%%%%%%%%%%%%%%%%%%%
\subsection{Theorem~\ref{thm:Matrix_Recovery_High_Prob}}
%%%%%%%%%%%%%%%%%%%%%%%%%%%%%%%%%%%%%%%%%%%%%%%%%%%%%%%%%%%%%%%%%%%%%%%%%%%%
\begin{IEEEproof}
Set $\bs{\Phi}=\bs{\Phi}_\Alltop$. Let $\mathscr{A}$ denote the
event that $\norm{\bs{X^*X} - \bs{I}} < \half$, and let
$\mathscr{B}$ represent the event that BP recovers random $\bs{s}$
from the observation $\bs{y} = \bs{H}\bs{f}_\Alltop =
\bs{\Phi}_\Alltop\bs{s}$. Proposition~\ref{prop:Tropp_Thm_B}
concerns $\Pr{\gr{\mathscr{A}^\complement}}$ where
$\mathscr{A}^\complement$ is the complement of set $\mathscr{A}$,
and Proposition~\ref{prop:Tropp_Thm_14} addresses
$\Pr{\gr{\mathscr{B}|\mathscr{A}}}$. To apply these propositions we
need their conditions to be satisfied simultaneously. Since
$\bs{\Phi}_\Alltop$ is a unit-norm tight frame we know that
$\norm{\bs{\Phi}_\Alltop}^2 = N$. With $\mu = {1}/{\sqrt{N}}$ and
taking $\delta = \half$ the condition of Proposition
\ref{prop:Tropp_Thm_B} is
\begin{equation} \label{eq:Thm_B_Alltop_condition}
\sqrt{\frac{K}{N}\log{K}\cdot\vartheta} \;+\; \frac{K}{N} \;\leq\;
\frac{c}{2}.
\end{equation}
Fix $\zeta = \varepsilon^2$ for some sufficiently small desired
probability of error in Proposition~\ref{prop:Tropp_Thm_14}. The
sparsity condition can now be rewritten as $K\leq
N/\gr{16\log{\gr{N/\varepsilon}}}$. Substituting this into
(\ref{eq:Thm_B_Alltop_condition}) the LHS is less than
\begin{eqnarray} \label{eq:Thm_B_Alltop_LHS_condition}
{}&& \!\!\!\!\!\!\!\!\!\!\!\!\!\!\!\!\!\!\!\!\!\!\!\!\!\!\!\!\!\!
\sqrt{\frac{\vartheta}{16\log{\gr{N/\varepsilon}}}
\log\Gr{\frac{N}{16\log{\gr{N/\varepsilon}}}}} \;+\; \frac{1}{16\log{\gr{N/\varepsilon}}} \nonumber \\
&<& \sqrt{\frac{\vartheta}{16\log{\gr{N/\varepsilon}}} \log{{N}}}
\;+\;
\sqrt{\frac{1}{16\log{\gr{N/\varepsilon}}}} \nonumber\\
&<&
\frac{1}{2}\sqrt{\frac{\vartheta\;\!{\log{{N}}}}{\log{\gr{N/\varepsilon}}}}
\qquad \mbox{(since $\vartheta, \log{N} \geq 1$)} .
\end{eqnarray}

Choose $\vartheta\geq1$ such that
$\sqrt{{\vartheta\;\!{\log{{N}}}}/{\log{\gr{N/\varepsilon}}}} \leq
c$ is satisfied. Assume the other conditions of
Proposition~\ref{prop:Tropp_Thm_14} (observe that event
$\mathscr{A}$ implies $\sigma_{\min}\gr{\bs{\Phi}_T}\geq
1/\sqrt{2}$), and let $\bs{X} = \bs{\Phi}_T$ in
Proposition~\ref{prop:Tropp_Thm_B}. Then
\begin{eqnarray} \label{eq:Thm_2_main_proof}
\Pr{\gr{\mathscr{B}}} &\geq& \Pr{\gr{\mathscr{B}|\mathscr{A}}} \;\! \Pr{\gr{\mathscr{A}}} \nonumber \\
&\geq& \gr{1-2\varepsilon^2} \gr{1-K^{-\vartheta}} \nonumber \\
&>& 1-2\varepsilon^2-K^{-\vartheta}
\end{eqnarray}
as desired.\end{IEEEproof}

%%%%%%%%%%%%%%%%%%%%%%%%%%%%%%%%%%%%%%%%%%%%%%%%%%%%%%%%%%%%%%%%%%%%%%%%%%%%
\subsection{Theorem~\ref{thm:Matrix_Recovery_Guaranteed_Noise}}
%%%%%%%%%%%%%%%%%%%%%%%%%%%%%%%%%%%%%%%%%%%%%%%%%%%%%%%%%%%%%%%%%%%%%%%%%%%%
\begin{IEEEproof}
As in the proof of Theorem~\ref{thm:Matrix_Recovery_Guaranteed},
this follows immediately \emph{mutatis mutandis}.
\end{IEEEproof}

%%%%%%%%%%%%%%%%%%%%%%%%%%%%%%%%%%%%%%%%%%%%%%%%%%%%%%%%%%%%%%%%%%%%%%%%%%%%
%%%%%%%%%%%%%%%%%%%%%%%%%%%%%%%%%%%%%%%%%%%%%%%%%%%%%%%%%%%%%%%%%%%%%%%%%%%%
\section*{Acknowledgment}
%%%%%%%%%%%%%%%%%%%%%%%%%%%%%%%%%%%%%%%%%%%%%%%%%%%%%%%%%%%%%%%%%%%%%%%%%%%%
%%%%%%%%%%%%%%%%%%%%%%%%%%%%%%%%%%%%%%%%%%%%%%%%%%%%%%%%%%%%%%%%%%%%%%%%%%%%
The authors would like to thank Roman Vershynin at UC Davis,
Benjamin Friedlander at UC Santa Cruz, Joel Tropp at the California
Institute of Technology, and Jared Tanner at the University of
Edinburgh for many fruitful discussions. Additionally, the authors
acknowledge and appreciate the useful comments and corrections from
the anonymous reviewers.

\ifCLASSOPTIONcaptionsoff
  \newpage
\fi

%%%%%%%%%%%%%%%%%%%%%%%%%%%%%%%%%%%%%%%%%%%%%%%%%%%%%%%%%%%%%%%%%%%%%%%%%%%%
%%%%%%%%%%%%%%%%%%%%%%%%%%%%%%%%%%%%%%%%%%%%%%%%%%%%%%%%%%%%%%%%%%%%%%%%%%%%
%%%%%%%%%%%%%%%%%%%%%%%%%%%%%%%%%%%%%%%%%%%%%%%%%%%%%%%%%%%%%%%%%%%%%%%%%%%%
\bibliography{../matt_refs}
\end{document}